 
\input amssym.def  
\input amssym.tex 
 
\def\dual{{^*}} 
\def\cG{{G}} 
\def\cL{{\cal L}} 
\def\cO{{\cal O}} 
\def\cT{{\cal T}} 
\def\cN{{\cal N}}

\def\bG{{\Bbb G }} 
\def\bQ{{\Bbb Q }} 
\def\bZ{{\Bbb Z }} 
 
\def\bC{{\Bbb C }} 
\def\bP{{\Bbb P }} 
\def\Hom{{\rm Hom}} 
\def\Slutt{\sqcap\!\!\!\!\sqcup} 
\def\qed{\hfill $\Slutt$ \medskip}

\font\tit=cmbx10 scaled\magstep2

{\centerline{ \tit Moduli of (1,7)-polarized abelian surfaces via syzygies}} 
 
\vskip10pt 
\centerline {\bf Dedicated to the memory of Alf B. Aure } 
\vskip20pt 
\centerline{\it By Nicolae Manolache\footnote *{partially supported by the 
Humboldt Gesellschaft and Max Planck Institute f\"ur Mathematik, Bonn} 
at Bucharest  and Frank-Olaf Schreyer at Bayreuth} 
\vskip20pt 
{\bf Abstract.} We prove that the moduli space $X(1,7)$ of 
(1,7)-polarized abelian surfaces with canonical level-structure is 
birational to the Fano 3-fold $V_{22}$ of polar hexagons of 
the Klein quartic $\overline X(7)$. In particular $X(1,7)$ is rational and 
the birational map to $\bP^3$ is defined over $\bQ$. As a byproduct we 
obtain explicitely the equations of the $(1,7)$--very-ample-polarized abelian 
surfaces embedded in $\bP^6$. 
  
 
\vskip20pt 
\hskip10pt\vbox{\settabs 2 \columns 
 \+ 0. Introduction & \cr 
 \+ 1. Review & \cr 
 \+ 2. Syzygies & \cr 
 \+ 3. Symmetry & \cr 
 \+ 4. Moduli & \cr 
 \+ A. Appendix & \cr 
 \+ References & \cr } 
\bigskip 
 
{\bf 0. Introduction.} 
 
Moduli spaces of polarized abelian varieties is a much studied subject. The  
common approach to their construction is as arithmetic quotient of the Siegel 
upper half space. Their study involves then the beautiful subject of modular 
forms. In this paper we follow a different approach. We construct them as Heisenberg 
invariant part of a Hilbert scheme. 
 
In the particular case of $(1,7)$-polarized abelian surfaces we can obtain in this way  
only a birational model of the moduli space, because not every polarization is very ample. 
However our method goes quite far: We obtain a birational parametrization of the 
moduli space defined over $\bQ$. 
 
Moreover we discovered a new relation between the moduli space $X(1,7)$ and the modular 
curve $X(7)$ of elliptic curves with a level 7 structure. The variety 
of sums of 6 powers of the Klein quartic $\overline X(7)$, ie. the
variety of polar hexagons of the Klein quartic, is our model of 
$X(1,7)$. 
 
\bigskip 
 
The paper is organised as follows, each section devoted to one basic idea. 
  In section 1 we review the construction of the moduli space as invariant part of the  
Hilbert scheme and collect some basic notations. 
 
In section 2 we recall that a very ample line bundle of class $(1,7)$ embeds 
an abelian surface projectively normal and study its syzygies. Due to  
$H^1(A,{\cal O}) \ne 0$ the minimal free resolution is longer than the codimension.  
However if we allow a locally free resolution, there is a rather natural self-dual  
resolution $F$. A result like this should hold quite generally for Gorenstein subvarieties  
of smooth manifolds, whose canonical bundle is induced, cf. [EPW], [W] 
for some results in this direction. 
 
Section 3 brings in the action of the Heisenberg group. With this, the middle 
syzygy map boils down to a $3 \times 2$ matrix, which can be interpreted as the  
Hilbert-Burch matrix (cf. 
[E] Thm 20.15) of a twisted cubic in a certain $\bP^3$. 
The complex condition on $F$ gives certain linear relations among the coefficients of  
the defining quadratic equations of this twisted cubic. This is enough to determine 
our model of the moduli space: It is a particular Fano 3-fold $V_{22} \subset \bP^{13}$ 
of degree 22. 
 
In section 4 we recall the various descriptions of a $V_{22}$, one of them being 
the variety of sums of powers of a plane quartic curve. For our $V_{22}$ this is  
the Klein quartic. We finish with the birational map $\bP^3 - \to X(1,7)$ induced 
by the triple projection from a particular point of the $V_{22}$. 
 
The Appendix contains some  formulas of the representation theory of the 
Heisenberg group $H_7$ and $SL_2(\bZ_7)$. 
 
\bigskip 
\noindent In several aspects we are not yet completely satisfied with our results here. 
 
\item{1.} A detailed study of the geometry of the surfaces (over boundary points) 
comparable to the study of   
Barth, Hulek and Moore [BHM] of $X(1,5)$ via the Horrocks-Mumford bundle 
 is missing. 
 
\item{2.} Points on our $V_{22}$ parametrize a family of 10-nodal Kummer quartics in  
$\bP^3$. Could it be that every $V_{22}$ parametrizes some 10-nodal quartics? 
 
\noindent However the paper had already some fruits: As observed first by Alf Aure and Kristian  
Ranestad, conics on the $V_{22}$ correspond to pencils of abelian surfaces which sweep 
out a Calabi-Yau 3-fold $Y$ of degree 14. In some sense the study of pairs 
$\{(A,Y) \mid A \subset Y \}$ of abelian surfaces contained in some Calabi-Yau is easier 
than studying $A$'s alone. This point of view turned out to be the key 
in the solution of Gross and Popescu [GP] of the problem posed 
by Gritsenko [Gr] to decide, which Siegel modular 3-folds are rational. 
\medskip 
{\bf Acknowledgement.} We thank Geir Ellingsrud for pointing us towards $V_{22}$'s, 
 Shigeru Mukai for sharing his insight into the $V_{22}$'s with 
 us, and Kristian Ranestad for various discussions about this material. 
 During various stages in the preparation of this paper  
N. Manolache was supported by a Humboldt fellowship at Bayreuth University 
and enjoyed the hospitality of Max Planck Institute in Bonn. He expresses his  
gratitude to these institutions. Special thanks are due to Professor  
F. Hirzebruch for the encouraging support.  
 
\medskip 
 
{\bf Notation.} 
Most of the notation will be introduced directly in the text. 
We recall here only some of it and also some notation which will be used 
tacitly: 
 
$A$ an abelian variety 
 
$\cL$ an ample line bundle of type $(1,7)$ on $A$ 
 
$\hat A=Pic^0(A)$ 
 
$t_x$ is the automorphism of translation by $x\in A$

$V\dual $ for the dual of a vector space $V$

$UV$ or $U\cdot V$ for $U\otimes V$, where $U$, $V$ are vector spaces 
 
$nV$ for $\oplus _1^n V$ 
 
$\bG(k, V) $ for the Grassmann variety of k dimensional subspaces of 
the vector space V 

$\bP(V)$ the projective space of lines in V.
 
\noindent We shall use the Macaulay short hand notation for numerical data of a 
free resolutions over the graded polynomial ring $R$ and for its 
sheafified version. A table like 
$$\matrix{  
1 & - & - & - & - \cr 
- & 7 & 8 & - & - \cr 
- & - & 3 & 8 & 3 \cr}$$ 
stands for a complex 
$$ R \leftarrow F_1 \leftarrow F_2   
\leftarrow F_3 \leftarrow 
F_4 \leftarrow 0$$ 
with 5 terms $F_0 = R$ and $F_i = \oplus_{k=1}^{r_i}R(-a_{ik})$ for 
$i=1,\ldots,4$ where the number 
of generators of $F_i$ in a given degree are encoded by the numbers 
in the $i^{th}$ column. More precisely the number in position $(i,-j)$ 
in the table is the number of generators of degree $i+j$ of $F_i$. 
In the example above, the image of $F_1=7R(-2)$ is an ideal generated by 
$7$ quadrics, which have $8$ linear syzygies and further $3$ 
quadratic syzygies, corresponding to $F_2 =8R(-3)\oplus 3R(-4)$; 
$F_3 = 8 R(-5)$, $F_4= 3R(-6)$. 
 
Notice that the entries of a block in a syzygy matrix above corresponding to 
two consecutive numbers in the same line are linear, while the maps 
to the upper left and from the lower right corners of a square are 
quadratic. 
Examples: The syzygies of the twisted cubic in $\bP^3$ have shape 
$$\matrix{  
1 & - & - \cr 
- & 3 & 2 \cr}\ \ \ \ .$$ 
The syzygies of a plane cubic in $\bP^3$ union a point, which is not in that 
plane, look like 
$$\matrix{  
1 & - & - & -\cr 
- & 3 & 3 & 1\cr 
- & 1 & 1 & -\cr}\ \ \ \ .$$

\bigskip 
{\bf 1. Review. } We review some well known facts about abelian varieties 
following [Mum66] and [Mum70] (see also [LB]).  
 
\medskip 
 
(1.1) Let $A$ be an abelian variety with an ample line bundle $\cL$, 
$$\phi_\cL : A \rightarrow Pic^0(A), \quad x \mapsto t_x^* \cL \otimes 
\cL^{-1}$$ the associated group homomorphism, and let $K_\cL$ be the kernel 
of $\Phi _\cL$.  $K_\cL$ is finite, since $\cL$ is ample (cf. [Mum70], 
II.6. Application 1).  Moreover, by the theorem of the cube, $\Phi 
_{t_y^*\cL} = \Phi _\cL$, i.e. $\Phi_\cL$ depends only on $c_1(\cL)$ (the {\it 
polarization} class of $\cL$).

By definition, $K_\cL$ operates on $\bP(H^0(A,\cL))$, but not on 
$H^0(A,\cL)$.  The pullback central extension group $\cG _\cL$ is called the 
(infinite) Heisenberg group of $\cL$: 
$$\matrix{ 1 & \to & \bC^* & \to & \cG _\cL & \to & K_\cL & \to & 0 \cr 
            &   &  &   &   &  &   &   & \cr 
            &   & \parallel &   & \downarrow &  & \downarrow &   & \cr 
           &   &  &   &   &  &   &   & \cr 
       1 & \to & \bC^* & \to & GL(H^0(A,\cL)) & \to &PGL(H^0(A,\cL)) & \to & 1  \cr}$$ 
 
We will work mainly with a finite version of $\cG _\cL$: Since $\cG _\cL$ is a 
central extension of an abelian group, taking commutators induces a 
skew bilinear map 
$$e^\cL \colon K_\cL \times K_\cL \to \bC^*, \quad e^\cL(x,y)=\tilde x \cdot \tilde y 
\cdot \tilde x^{-1} \cdot \tilde y^{-1}\ \quad (\tilde x,\tilde y \in \cG_\cL \quad 
\hbox{preimages of}\ x,y\ ).$$ 
 
Since $K_\cL$ is finite, we may replace $\bC^*$ by the finite group 
$\mu_\cL$ generated by the image of $e^\cL$. Moreover $\mid \mu_\cL \mid$ is 
a divisor of the exponent of $K_\cL$.  The (infinite) Heisenberg group 
$\cG _\cL$ is the pushout of a certain (finite) Heisenberg group $H_\cL$ 
: 
 
$$\matrix{ 1 & \to & \mu_\cL & \to & H_\cL & \to & K_\cL & \to & 0 \cr  
           &   &  &   &   &  &   &   & \cr 
        &   & \downarrow &   & \downarrow &  & \downarrow &   & \cr 
           &   &  &   &   &  &   &   & \cr 
       1 & \to & \bC^* & \to & \cG _\cL & \to & K_\cL & \to & 0 \quad .\cr}$$

\medskip 
 
(1.2) An explicit description of $H_\cL$ uses more notation. In the 
analytic context $c_1(\cL) \in H^2(A^{an}, \bZ) \cong \Lambda ^2 
Hom(\Gamma, \bZ)$, where $ \Gamma$ is a lattice in $\bC ^g$ such that 
$A \cong \bC ^g/ \Gamma$. Hence $c_1(\cL)$ is an alternating 
nondegenerated 2-form on $ \Gamma$ with values in $\bZ$, (cf. [Mum70] 
p. 16).  In a convenient basis of $ \Gamma$ this alternating form can 
be written as a matrix: 
$$ c_1(\cL) =  \pmatrix {0 & \Delta \cr 
            - \Delta & 0\cr }  
$$ 
where $\Delta$ is a diagonal matrix $(d_1,\ldots,d_g)$ of positive 
integers such that $d_1 \mid d_2 \mid \ldots \mid d_g $. The 
collection of elementary divisors $(d_1,\ldots , d_g) =:\delta $ is 
called the type of the line bundle $\cL$. 
 
The diagram 
$$\matrix{ 0 &\to &\Gamma &\to &\bC^g & \to & A & \to & 0 \cr  
             &    &       &    &      &         &     &   & \cr 
             &c_1(\cL)&\downarrow &  &\downarrow \cong &  &\downarrow &\Phi_\cL & \cr 
             &   &  &   &   &  &   &   & \cr 
           0 &\to &\hat \Gamma & \to & \hat \bC^g & \to & \hat A & \to & 0 \quad \cr}$$ 
 
\noindent yields an isomorphism $K_\cL \cong {\hat \Gamma / \Gamma} \cong (\bZ/d_1 
\bZ)^{ 2} \oplus \ldots \oplus (\bZ/d_g \bZ)^{ 2}$. 
 
In the algebraic setting the elementary divisors of $K_\cL$ define the 
type: $e^\cL$ is nondegenerated and produces a decomposition of $K_\cL$ as 
a direct sum of two subgroups $K_1$, $K_2$, where $K_2 \cong \hat K_1 
\cong Hom(K_1,\bC^*)$, (cf. [Mum66], p. 293).  The most convenient 
description of $H_\cL$, up to isomorphism, is as the following extension 
$H_{\delta}$: 
$$  
1 \to \mu_{d_g} \to H_{\delta} \to ((\bZ/d_1 \bZ) \oplus \ldots \oplus ((\bZ/d_g  
\bZ)) \oplus (\mu_{d_1} \times \ldots \times \mu_{d_g}) \to 0 \quad , 
$$ 
where $\mu _{d} = Hom(\bZ/d\bZ,\bC^*)$ is the group of the $d$-th 
roots of unity. If we write the last term in the above exact sequence 
as $K(\delta ) \oplus   \hat{K}(\delta )$, with $\hat K(\delta 
)=Hom(K(\delta ), \bC ^*)$, then the multiplication in $H_{\delta 
}=\mu _{d_g} \oplus K(\delta) \oplus   \hat{K}(\delta )$ is defined by 
$$  
(\alpha , x, \rho) \circ (\beta ,y,\sigma) = (\alpha \cdot \beta \cdot  
\sigma (x),x+y,\rho \cdot \sigma). 
$$ 
On $K(\delta ) \oplus  \hat{K}(\delta )$ one has also naturally an alternate 
multiplicative form $e^\delta $. If $e_1,\dots ,e_g$ and $e_1^\ast ,\dots ,e_g^\ast$ 
are the canonical basis of $K(\delta )$ and $\hat K(\delta )$ respectively  
then  
$$ 
e^\delta : (K(\delta ) \oplus  \hat{K}(\delta ) \,)\times   
(K(\delta ) \oplus   \hat{K}(\delta ) \,)\to \mu _{d_g}  
$$ 
is defined by: 
$$  
e^\delta (e_r^\ast ,e_r) = (e^\delta (e_r,e_r^\ast))^{-1}=\exp ({2\pi i}/d_r), \ \ 
\hbox{for all}\ r=1,\ldots ,g \  
\hbox{and}\ 1\ \hbox{otherwise}  
$$ 
 
\medskip 
 
(1.3) $H^0(A,\cL)$ is an irreducible $H_\cL$-module. The argument uses 
that $K(\delta)$ lifts to a subgroup of $H_\cL$, and that $\cL$ descends 
to $A/K(\delta)$, (cf. [Mum66], pp. 290, 297). $H^0(A,\cL)$ is the unique 
irreducible representation of $H_\cL$ , on which the center $\mu_\cL 
\subset \bC^*$ acts by scalar multiplication. This $(\prod_{i=1}^{g} 
d_i)$-dimensional representation V is called the Schr\"odinger 
representation of $H_{\delta}$. 
 
\medskip 
  
(1.4) {\bf Definition.} ([Mum66]) {\it A {\bf theta--structure} 
on the pair $(A,\cL)$ is any isomorphism $\alpha:~H_\cL~\to~H_{\delta} $ which 
induces the identity on the centers $\mu_\cL$ and $\mu_{d_g}$ viewed as 
subgroups of $\bC^*$.  A theta--structure induces a {\bf level-structure  
of canonical type} i.e. 
a symplectic isomorphism $\alpha': K_\cL \to K_{\delta}$, symplectic with 
respect to $e^\cL$ and $e^\delta$. 
 
\noindent Since $\Phi_\cL$ depends only on the numerical equivalence class of  
$\cL$ we can speak of a level-structure for a polarized abelian variety.} 
 
The reason to consider level-structures is the following: 
\medskip 
 
(1.5) {\bf Theorem.}(Mumford) {\it There exists a coarse moduli space  
${\cal M}(\delta)$ for ample  
polarized abelian varieties with level-structure of type 
 $\delta =(d_1,\ldots,d_g)$.} 
 
\medskip 
 
(1.6) If $d_1 \ge 3$ the proof of this is elementary: By Lefschetz's 
theorem an ample line bundle $\cL$ of type $\delta =(d_1,\ldots,d_g)$ 
with $d_1 \ge 3$ is very ample.  The irreducibility of the 
Schr\"odinger representation $V$ and the level-structure on $(A,\cL)$ 
gives a canonical identification $\bP(V) = \bP(H^0(A,\cL))$ by Schurs 
lemma.  Thus every pair (A,\cL) with level structure of type $\delta$ 
occurs as a point in the Hilbert scheme $Hilb(\bP(V))$. Moreover $A 
\subset \bP(V)$ is $ H_\delta$-invariant. So we can take as ${\cal 
M}({\delta})$ an open part of an irreducible component of the fixpoint 
set  $ Hilb(\bP(V))^{H_{\delta}}$. To see that this has the right 
(in fact reduced) scheme structure we compare the tangent spaces: 
$$T_{Hilb^{H_\delta},A} = H^0(A,\cN_A)^{H_{\delta}} \cong  
Im(H^0(A,\cN_A) \to H^1(A,\cT_A) ).$$ 
Indeed, from the exact sequences: 
$$ 0 \to \cT_A \to \cT \to \cN_A \to 0 $$ 
$$0 \to \cO_A \to V\cO_A(1) \to \cT \to 0$$ 
where $\cal T$ is the restriction to $A$ of the tangent bundle of 
$\bP(V)$, one deduces the exact sequence: 
$$ 
\matrix{  
0 & \to & H^0(\cT_A) & \to & H^0(\cT) & \to & H^0(\cN_A) & \to H^1(\cT_A) 
& \to & \ldots \cr 
& & & & & & & & & & \cr 
& & \parallel & & \parallel & & &\parallel & & & \cr 
& & & & & & & & & & \cr 
& &  gI &  & gI\oplus Z & & &{g(g-1)\over 2}I & & &  \cr} 
$$ 
where $I$ is the trivial $1$--representation of $H_\delta$ and $Z$ is the complement 
of $gI$ in $H^0(\cT)$. 
This shows that that the image of $H^0(\cN_A)=Z\oplus H^0(\cN_A)^{H_\delta}$ 
in $H^1(\cT_A)$ is $H^0(\cN_A)^{H_\delta}$. \qed

(1.7) {\bf Remarks.} (1) Notice that the universal family over the 
Hilbert scheme is not necessarily a universal family in the sense of moduli, 
because, as subscheme, $A$ has no distinguished origin. Indeed picking 
the origin appropriate we obtain that $\cO(1)$ restricts to any 
desired line bundle $\cL$ within the polarization class $c_1(\cL)$.

(2) In general $Hilb(\bP(V))^{H_{\delta}}$ has many components. For 
example in the elliptic curve case, only if $d_1$ is prime, there is a 
single component.  For composed numbers, there are several components 
whose points correspond to union of $d_1/n$ elliptic curves of degree 
n, for every divisor n of $d_1$, (cf. [P]).

\bigskip 
 
{\bf 2. Syzygies.}  
 
\medskip 
 
(2.1) According to a theorem of Reider, a line bundle of type 
$(1,d),\ d\ge 5$ on a general abelian surface  
is  very ample (cf. [R] or [LB] p. 301-302). In this section we describe the  
syzygies of an abelian surface $A \subset \bP^6$ embedded by an very ample line  
bundle of type (1,7). Although the canonical class $\omega_A$ is induced from  
$\bP^6$ the minimal free resolution is not symmetric, since the coordinate ring is 
not projectively Cohen-Macaulay. However the only obstacle for this is 
$H^1(A,\cal O)$.  Using also some locally (not globally) free  sheaves in the 
resolution we obtain a nice self-dual resolution.

\medskip 
(2.2) Let $(A,\cL)$ be an abelian surface with a very ample line bundle 
of type $(1,7)$, and consider its image $$ A \hookrightarrow \bP^6 = 
\bP(V),$$ as an $H_7 = H_{(1,7)}$-invariant subvariety. 
 
\medskip 
(2.3) {\bf Lemma.} {\it $A \subset \bP^6$ is not contained in a quadric.} 
\medskip 
 
{\it Proof.} If a quadric would contain $A$ then $h^0(\bP^6,{\cal 
I}_A(2)) \ge 7$, since $H^0(\bP^6,\cO(2))$ is a direct sum of 
irreducible representations of dimension $7$, cf. Appendix. This is 
too much: By Castelnuovo's argument the 14 points $Z = A \cap \bP^4 $ 
for a general $\bP^4 \subset \bP^6$ impose at least 9 conditions on 
quadrics, i.e. $h^ 0(\bP ^4,{\cal I} _ Z (2)) \le 6$. On the other hand 
$h^0(\bP^4,{\cal I}_Z(2))=h^0(\bP^6,{\cal I}_A(2))+2 \ge 9$ by the 
exact sequence 
$$ 0 \rightarrow {\cal I}_A \rightarrow 2{\cal I}_A(1) \rightarrow  
{\cal I}_A(2) \rightarrow {\cal I}_{Z,\bP^4}(2) 
\rightarrow 0, \leqno{(2.3.1)}$$ 
a contradiction. \qed

(2.4) {\bf Corollary.} {\it $A \subset \bP^6$ has syzygies: 
$$\matrix{  
1 & - & - & - & - & - \cr 
- & - & - & - & - & - \cr 
- & 21 & 49 & 42 & 14 & 2 \cr 
- & - & - & - & 1 & - &  \cr}$$} 
\medskip
           
{\it Proof.} The map $H^0(\bP^6,\cO(2)) \to H^0(A,\cL^{\otimes2})$ 
is injective by the Lemma, hence an isomorphism, because both spaces 
are 28-dimensional.  So $H^1(\bP^6,{\cal I}_{A}(2)) = 0 $, i.e. $A 
\subset \bP^6$ is quadratically normal.  By the above: 
$h^0(\bP^4,{\cal I}_Z(2))=2$. So $Z \subset \bP^4$ has the Hilbert 
function $$(1,5,13,14,14,\ldots),$$ because its different function has 
no negative value, being the Hilbert function of an artinian ring. It 
follows $H^1(\bP^4,{\cal I}_{Z}(n)) = 0 $ for $n \ge 3$ and induction 
with the sequence (2.3.1) gives: $A \subset \bP^6$ is projectively 
normal. 
 
Moreover $\cal I_A$ is 4-regular in the sense of Castelnuovo-Mumford 
and nonzero syzygy numbers can only be in the indicated range, some of 
whose values are clear: 
$$\matrix{  
1 & - & - & - & - & - \cr 
- & - & - & - & - & - \cr 
- & 21 & ? & ? & ? & 2 \cr 
- &  ? &  ? &  \hbox{\rm \char'076} & 1 & - & \cr}$$ 
Namely, the number of cubic generators of the ideal is 
$h^0(\bP^6,{\cal I}_{A}(3)) = h^0(\bP^6,\cO(3))-h^0(A,{\cal 
O}_{A}(3)) =21 $.  The last $2$ represents $h^1(\cO_A)$ and the 
$1$ comes from the facts that $\omega _A \cong \cO_A$ and that
dualizing the above resolution one 
obtains ${\cal E}xt^4_{\cO_{\bf P}}(\cO_A,\omega _{\bf P}) 
\cong \omega _A$. 
 
The last argument gives also 
$\hbox{\rm \char'076} = 0$, since  ${\cal E}xt^i_{\cO_{\bf P}}(\cO_A,\omega  
_{\bf P}) = 0$ for $i \le 3$  
and A is non-degenerate. Now all the vanishing is clear, and  
the nonzero values can be computed from the Hilbert function.  
\qed

The resolution has $length > codim A$, since $A \subset \bP^6$ is 
not arithmetically Cohen-Macaulay. In particular it is not 
symmetric. However, if we allow locally free sheaves instead of only direct 
sums of line bundles, then there is a nice self-dual resolution: 
\medskip 
 
(2.5) {\bf Theorem.} {\it $A \subset \bP^6$ has a self-dual resolution 
of type $$  
0\longleftarrow \cO_A \longleftarrow \cO 
{\buildrel \beta \over \longleftarrow} 21\cO(-3) {\buildrel \alpha 
\over\longleftarrow} 2\Omega ^3 {\buildrel \alpha' \over 
\longleftarrow }21\cO(-4) 
{\buildrel \beta'\over\longleftarrow}\cO(-7) \longleftarrow 0  
$$  
with $\alpha'=\pmatrix{ 0  & 1 \cr 
             -1  & 0}\ ^t\alpha$  
and $\beta'=\ ^t\beta$. 
} 
\medskip
 
{\it Proof.} The above resolution is obtained by a kind of 
subtracting a piece of the Koszul sequence multiplied with 
$h^1(A,\cO_A) $ from the resolution in Corollary (2.4). Let $\cal 
K$ be the kernel of the map $\cO \longleftarrow 21\cO(-3)$ 
in the resolution in (2.4). Then comparing the Koszul resolution of 
the $\bC$ vector space $H^1(A,\cO_A)\dual \cong Ext^5_S(S_A,S(-7))$ 
with the dual complex of (2.4), yields a commutative diagram with 
exact rows and columns: 
 
\vskip10pt 
\vbox  
{\tabskip =-1pt \offinterlineskip 
\halign to 400pt 
{ 
\strut # & \hfil # \hfil & \hfil # \hfil & \hfil # \hfil & \hfil # \hfil & \hfil  
# \hfil & \hfil # \hfil & \hfil # \hfil & \hfil # \hfil & \hfil # \hfil & \hfil #  
\hfil & \hfil # \hfil & \hfil # \hfil \cr 
   &   &    &    &        0          &    &     &     &              &     &      
   &    &  \cr 
   &   &    &    &                   &    &     &     &              &     &      
   &    &  \cr 
   &   &    &    & $\downarrow $     &    &     &     &              &     &      
   &    &  \cr 
   &   &    &    &                   &    &     &     &              &     &     
    &    &  \cr 
   &   &    &    &$21\cO(-4)$ &    &     &     &     $0$      &     &     &     
   &  \cr 
   &   &    &    &                   &    &     &     &              &     &      
   &    &  \cr 
   &   &    &    &  $\downarrow $    &    &     &     &$\downarrow$  &     &      
   &    &  \cr 
   &   &    &    &                   &    &     &     &              &     &      
   &    &  \cr 
0  &$\leftarrow$ &  $2 \cdot \Omega ^3$ & $\leftarrow $ &  
  $2\cdot 35\cO(-4)$ &  $\leftarrow $ &  $2\cdot 21\cO(-5)$  &  $\leftarrow $ &  
  $2\cdot 7\cO(-6) $ &  $\leftarrow $ &  $2\cdot \cO(-7) $ &   
 $\leftarrow $ &  0 \cr 
&   &   &   &   &   &   &   &   &   &   &   &  \cr 
&   &   &   & $\downarrow $ &  &  $\downarrow \!{\cal o}\  $ &   & $\downarrow $ &   
& $\parallel $ &   &  \cr 
&   &   &   &   &   &   &   &   &   &   &   &  \cr 
 0  &  $\leftarrow $ &  $\cal K$ & $\leftarrow $ &  
 $49\cO(-4) $ &  $\leftarrow $ &  $42\cO(-5) $ &   
$\leftarrow $ & $14\cO(-6) \oplus \cO(-7)$  & $\leftarrow$ & $2\cO(-7)$ &   
$\leftarrow $ & $0$ \cr 
   &   &    &    &                   &    &     &     &               &     &     &   
     &  \cr 
   &   &    &    &  $\downarrow$     &    &     &     &$\downarrow$   &     &     &   
     &  \cr 
   &   &    &    &                   &    &     &     &               &     &     &   
     &  \cr 
   &   &    &    &        $0$        &    &     &     &$\cO(-7)$ &     &     &    &   
   \cr 
   &   &    &    &                   &    &     &     &               &     &     &  
      &  \cr 
   &   &    &    &                   &    &     &     &$\downarrow$   &     &     &   
     &  \cr 
   &   &    &    &                   &    &     &     &               &     &     &  
      &  \cr 
   &   &    &    &                   &    &     &     &      $0$      &     &     &   
     &  \cr 
\cr} 
} 
 
The map $2 \cdot 21\cO(-5) \to 42\cO(-5)$ is surjective, and 
$ ker(2 \cdot 35\cO(-4) \to 49\cO(-4)) \cong 21{\cal 
O}(-4)$, because otherwise $A \subset \bP^6$ would be contained in a 
quadric, or more then 21 cubics. A diagram chase gives the desired 
resolution.  To see the assertions about the maps, we compare this 
complex with its dual: 
\bigskip 
\hskip5pt 
\vbox  
{
\halign to 420pt 
{\strut # & \hfil # \hfil & \hfil # \hfil & \hfil # \hfil & \hfil # \hfil & \hfil # \hfil & \hfil # \hfil & \hfil # \hfil & \hfil # \hfil & \hfil # \hfil & \hfil # \hfil & \hfil # \hfil & \hfil # \hfil & \hfil # \hfil & \hfil # \hfil \cr 
  0  & $\leftarrow$ &  $\cO_A$ & $\leftarrow $ &  
  $\cal O$ &  ${\buildrel \beta \over \leftarrow} $ &  $21\cO(-3)$  &   
${\buildrel \alpha \over \leftarrow} $ &  
$H^1(\cO_A)\otimes\Omega^3 $ &  ${\buildrel \alpha'\over \leftarrow} $ &   
$21\cO(-3) $ & ${\buildrel \beta' \over \leftarrow} $ &  $\cO(-7)$ &  
$\leftarrow $& $0$ \cr 
&   &   &   &   &   &   &   &   &   &   &   &  &   & \cr 
&   & $\downarrow \!{\cal o}\ \phi $ &   & $\downarrow \!{\cal o} $  &   & $ 
\downarrow \!{\cal o} $  &   & $\downarrow \!{\cal o}\ u $ & 
& $\downarrow \!{\cal o} $    &   & $\downarrow \!{\cal o} $   &  &    \cr 
&   &   &   &   &   &   &   &   &   &   &   &  &   & \cr 
0  & $\leftarrow$ &  $\omega _A$ & $\leftarrow $ &  
  $\cal O$ &  ${\buildrel ^t\beta' \over \leftarrow} $ &  $21\cO(-3)$  &  
${\buildrel ^t\alpha' \over \leftarrow} $ &  
$H^1(\cO_A)\dual\otimes\Omega^3 $ &  ${\buildrel ^t\alpha \over \leftarrow} $ &   
$21\cO(-3) $ & ${\buildrel ^t\beta \over \leftarrow} $ &  $\cO(-7)$ &  
$\leftarrow $& $0$ \cr 
\cr} 
} 
 
\bigskip 
 
The isomorphism $u$ is compatible with Serre duality. The diagrams: 
$$ 
\matrix{ 
H^1(\cO_A) &  \cong  & H^1(\cO_A) \otimes \bC \cr 
& & \cr 
H^1(\phi)\ \downarrow \!{\cal o}\  & & \downarrow \!{\cal o}\ \ u \cr 
& & \cr 
H^1(\omega_A) & {\buildrel \cong \over \longrightarrow} &   H^1(\cO_A)\dual \otimes \bC 
}  
\quad \hbox{and} \quad 
\matrix{ 
H^1(\cO_A)\otimes H^1(\cO_A) & \longrightarrow & H^2(\cO_A) \cr 
& & \cr 
\downarrow id\otimes H^1(\phi ) & & \downarrow  H^2(\phi ) \cr 
& & \cr 
H^1(\cO_A)\otimes H^1(\omega _A) & \longrightarrow &  H^2(\omega _A) 
} 
$$ 
commute.  So the map in the first row of the last diagram is antisymmetric, as the one in 
the second row is.  This gives the relations between $\alpha'$, $\beta'$ and 
$\alpha$, respectively $\beta$. 
\qed 
 
\bigskip 
 
{\bf 3. Symmetry.} Taking into account the symmetries of $(A,\cL)$ we 
have the following $H_7$ or $G_7$-invariant resolutions. 
\medskip 
 
(3.1) {\bf Corollary.} {\it Taking the canonical $H_7$-invariant 
embedding of $A$ corresponding to the taken polarization and the level 
structure, one gets: $$ 
0\leftarrow {\cal I}_A  \leftarrow 3V_4\cO(-3) \leftarrow 7V_1\cO(-4) 
\leftarrow 6V_2\cO(-5) \leftarrow 2V\cO(-6)\oplus \cO(-7) 
\leftarrow 2\cO(-7) \leftarrow 0  
$$ 
If one considers an $G_7$-invariant embedding one obtains: 
$$ 
0\leftarrow {\cal I}_A  \leftarrow 3V_4\cO(-3) \leftarrow (5V_1 \oplus 
2V_1^\sharp )\cO(-4) 
\leftarrow 6V_2^\sharp \cO(-5) \leftarrow 2V^\sharp \cO(-6)\oplus \cO(-7) 
\leftarrow 2S\cO(-7) \leftarrow 0  
$$ 
} 
 
{\it Proof.} Everything follows using the tables from the appendix. 
\qed

(3.2) {\bf Theorem.} {\it An abelian surfaces $G_7$-invariantly embedded 
in $\bP^6$ has an $G_7$-invariant resolution of the form: $$ 
0\longleftarrow \cO_A  \longleftarrow \cO {\buildrel \beta \over \longleftarrow}  
3V_4\cO(-3) {\buildrel \alpha \over\longleftarrow} 2S\Omega ^3  {\buildrel \alpha'  
\over \longleftarrow }3V_1\cO(-4) 
{\buildrel \beta'\over\longleftarrow}\cO(-7) \longleftarrow 0  
$$ 
with 
$\alpha'=\pmatrix{ 0  & 1 \cr 
             -1  & 0}\ ^t \alpha$  
and $\beta'=\ ^t \beta$. 
} 
 
\medskip 
 
{\it Proof.} The big diagram in the 
proof of theorem 2.5 reads as $G_7$-modules 
\bigskip  
 
\vbox  
{ 
\tabskip =-1pt \offinterlineskip 
\halign to 400pt 
{ 
\strut # & \hfil # \hfil & \hfil # \hfil & \hfil # \hfil & \hfil # \hfil & \hfil #  
\hfil & \hfil # \hfil & \hfil # \hfil & \hfil # \hfil & \hfil # \hfil & \hfil # \hfil  
& \hfil # \hfil & \hfil # \hfil \cr 
   &   &    &    &        0          &    &     &     &              &     &     &   
     &  \cr 
   &   &    &    &                   &    &     &     &              &     &     &   
     &  \cr 
   &   &    &    & $\downarrow $     &    &     &     &              &     &     &   
     &  \cr 
   &   &    &    &                   &    &     &     &              &     &     &    
    &  \cr 
   &   &    &    &$3V_1\cO(-4)$ &    &     &     &     $0$      &     &     &    &   
   \cr 
   &   &    &    &                   &    &     &     &              &     &     
    &    &  \cr 
   &   &    &    &  $\downarrow $    &    &     &     &$\downarrow$  &     &     &    
    &  \cr 
   &   &    &    &                   &    &     &     &              &     &     &   
     &  \cr 
0  &$\leftarrow$ &  $2S\cdot \Omega ^3$ & $\leftarrow $ &  
  $2S\cdot (4V_1^\sharp \oplus V_1)\cO(-4)$ &  $\leftarrow $ &  $2S\cdot (3V_2)\cO(-5) 
  $  &  $\leftarrow $ & $2S\cdot V\cO(-6) $ &  $\leftarrow $ &  $2S\cdot \cO(-7) $ &   
 $\leftarrow $ &  0 \cr 
&   &   &   &   &   &   &   &   &   &   &   &  \cr 
&   &   &   & $\downarrow $ &  &  $\parallel $ &   & $\downarrow $ &  & $\parallel $  
&   &  \cr 
&   &   &   &   &   &   &   &   &   &   &   &  \cr 
 0  &  $\leftarrow $ &  $\cal K$ & $\leftarrow $ &  
 $(5V_1\oplus 2V_1^\sharp)\cO(-4) $ &  $\leftarrow $ &  $6V_2^\sharp \cO(-5) $ &   
$\leftarrow $ & $2V^\sharp\cO(-6) \oplus \cO(-7)$  & $\leftarrow$ & $2S\cdot \cO(-7)$ 
 &  $\leftarrow $ & $0$ \cr 
   &   &    &    &                   &    &     &     &               &     &     &   
     &  \cr 
   &   &    &    &  $\downarrow$     &    &     &     &$\downarrow$   &     &     &   
     &  \cr 
   &   &    &    &                   &    &     &     &               &     &     &   
     &  \cr 
   &   &    &    &        $0$        &    &     &     &$\cO(-7)$ &     &     &    & \cr 
   &   &    &    &                   &    &     &     &               &     &     &   
     &  \cr 
   &   &    &    &                   &    &     &     &$\downarrow$   &     &     &   
     &  \cr 
   &   &    &    &                   &    &     &     &               &     &     &    
   &  \cr 
   &   &    &    &                   &    &     &     &      $0$      &     &     &  
     &  \cr 
\cr} 
} 
\bigskip 
and the result is clear. 
\qed

 Next we view $\alpha$ as a $3\times2$-matrix with entries in $\Hom(S\Omega 
 ^3, V_4\cO(-3))$. Since $\alpha$ defines a $G_7$-morphism the entries 
 lie in the $G_7$-invariant part. 
\medskip 
 
(3.3) {\bf Proposition.}{\it  $$\Hom _{G_7}(S\Omega ^3,V_4\cO(-3))=4I,$$ 
i.e. $\alpha$ has entries in a 4-dimensional vector space.} 
\medskip 
 
{\it Proof.}  $\Hom (\Omega ^3, \cO(-3))\cong \Lambda ^3 V 
=V_1\oplus 4V_1^\sharp$. Hence $ \Hom(S\Omega ^3, V_4{\cal 
O}(-3))\cong V_4\otimes (V_1^\sharp \oplus 4V_1) = 4I\oplus S \oplus 
5Z$ and $\Hom _{G_7}(S\Omega ^3,V_4\cO(-3))=4I$. 
 \qed

(3.4) {\bf Remark.} If ${\cal F}_1$ and ${\cal F}_2$ are two 
$G_7$-sheaves, $\Hom _{G_7}({\cal F}_1, {\cal F}_2)$ is a $N$-module, 
because $G_7=H_7\rtimes \bZ _2$ is a normal subgroup of $N \cong H_7\rtimes 
SL_2(\bZ_7)$, $\iota $ being central in $SL_2(\bZ_7)$. 
 
Using the character table of $SL_2(\bZ_7)$, one sees that 
$\Hom_{G_7}(S\Omega ^3, V_4\cO(-3)) \cong U'$, with the notation from 
the appendix. 
\medskip 
For the following considerations we choose a basis $u_0, \ldots , u_3$ 
of $U'$, so that in the decomposition into irreducible $G_7$-modules of 
$\Lambda ^3 V = V_1 \oplus (U'\otimes V_1 ) $, $u_0, \ldots , u_3 $ 
correspond to the $V_1 $ pieces generated as a $H_7$-module by 
$e_1\wedge e_4 \wedge e_2 - e_6\wedge e_3 \wedge e_5$, $e_0\wedge e_1 
\wedge e_6 $, $e_0\wedge e_2 \wedge e_5 $, or $e_0\wedge e_4 \wedge 
e_3 $ respectively.  Then, in the above decomposition of $\Lambda 
^3V$, $V_1$ is generated as a $H_7$-module by $e_1\wedge e_4 \wedge e_2 
+ e_6\wedge e_3 \wedge e_5$. More precisely, the above elements 
correspond to $ u_0\otimes e_0, \ldots , u_3\otimes e_0 $ in $U'\otimes 
V_1$ and the elements $u_k \otimes e_\ell $ are obtained permuting the 
indices of $e$'s via $\sigma $. 
 
\medskip

With these notations the matrix $\alpha =(a_{ij })$ will have entries $ 
(a_{ij }) =\sum_{k=0}^3 a_{ij }^k u_k$. We want to express more 
conveniently the condition $\alpha \alpha'=0$, where  
$\alpha'= \pmatrix{0 & 1 \cr -1 & 0 } \alpha$. 
 
\medskip 
 
(3.5) {\bf Proposition. } {\it A matrix $\alpha$ as above satisfies $\alpha \alpha'=0 
$ iff the three quadrics in $\bP^3 =\bP(U)$ given by its 
$2\times2$-minors are annulated by each of the three operators: 
$$ 
 \Delta _1 ={\partial ^2 \over \partial u_0\partial u_1} - {1 \over 2} {\partial ^2  
 \over \partial u_2^2}, \quad \quad 
\Delta _2 ={\partial ^2 \over \partial u_0\partial u_2} - {1 \over 2} {\partial ^2  
\over \partial u_3^2}, \quad \quad 
 \Delta _3 ={\partial ^2 \over \partial u_0\partial u_3} -  {1 \over 2}   
 {\partial ^2 \over \partial u_1^2}.  
$$ 
} 
 
{\it Proof. } We use the fact that under the identifications $\Hom 
(\cO (-4), \Omega ^3) =\Lambda ^3 V$, $\Hom (\Omega ^3, \cO 
(-3) ) =\Lambda ^3 V$ the composition of two maps $\cO(-4)\to 
\Omega ^3$, $ \Omega ^3 \to \cO(-3) $ is given by wedge product, 
if we identify canonically $\wedge ^6 V$ with $V\dual = V_3 = H^0({\cal 
O}(1))$. 
 
Observe now that $u_0$ interpreted as an element in $\Hom (V_1\cO 
(-4), S\Omega ^3)$ is given by the following $1\times 7$ matrix with 
entries in $\Lambda ^3 V$: 
 
$$ 
u_0=(e_{1+k}\wedge e_{4+k} \wedge e_{2+k}-e_{6+k} \wedge e_{3+k} \wedge  
e_{5+k} )_{k \in \bZ_7}  
$$ 
and similarly: 
 
$$ 
u_1=(e_k\wedge e_{1+k} \wedge e_{6+k})_k \ , \quad  
u_2=(e_k\wedge e_{2+k} \wedge e_{5+k})_k \ , \quad 
u_3=(e_k\wedge e_{4+k} \wedge e_{3+k})_k \ . 
$$ 
  
The same elements, interpreted in $Hom (S\Omega ^3 , V_4{\cal 
O}(-3))$, will be identified with the transpose of the above ones. 
Then the only compositions of two $u_i$'s which are not $0$ are: 
 
$$ 
u_0u_1 =u_1u_0 = -u_2u_2 = B_1 := \pmatrix{ 
0 & x_4 & 0 & 0 & 0 & 0 & -x_3 \cr 
-x_4 & 0 & x_5 & 0 & 0 & 0 & 0 \cr 
0 & -x_5 & 0 & x_6 & 0 & 0 & 0 \cr 
0 & 0 & -x_6 & 0 & x_0 & 0 & 0  \cr 
0 & 0 & 0 & -x_0 & 0 & x_1 & 0  \cr 
0 & 0 & 0 & 0 & -x_1 & 0 & x_2 \cr 
x_3 & 0 & 0 & 0 & 0 & -x_2 & 0} 
$$ 
 
\bigskip 
 
$$ 
u_0u_2 =u_2u_0 = - u_3u_3 = B_2 := \pmatrix{ 
0 & 0 & x_1 & 0 & 0  & -x_6 & 0 \cr 
0 & 0 & 0 & x_2 & 0 & 0 & -x_0 \cr 
-x_1 & 0 & 0 & 0 & x_3 & 0 & 0  \cr 
0 & -x_2 & 0 & 0 & 0 & x_4 & 0   \cr 
0 & 0 & -x_3 & 0 & 0 & 0 & x_5   \cr 
x_6 & 0 & 0 & -x_4 & 0 & 0 & 0 \cr 
0 & x_0 & 0 & 0 & -x_5 & 0 & 0 } 
$$ 
 
\bigskip 
 
$$ 
u_0u_3 =u_3u_0 = -u_1u_1 =B_3 := \pmatrix{ 
0 & 0 & 0 & -x_5 & x_2 & 0 & 0 \cr 
0 & 0 & 0 & 0 & -x_6 & x_3 & 0 \cr 
0 & 0 & 0 & 0 & 0 & -x_0 & x_4\cr 
x_5 & 0 & 0 & 0 & 0 &  0 & -x_1 \cr 
-x_2 & x_6 & 0 & 0 & 0 & 0 & 0  \cr 
0 & -x_3 & x_0 & 0 & 0 & 0 & 0 \cr 
0 & 0 & -x_4 & x_1 & 0 & 0 & 0} . 
$$ 
 
In particular $u_iu_j=u_ju_i$, i.e. the maps commute. Now, in the 
interpretation of $\alpha = (a_{ij})$ as a $3 \times 2$ matrix 
each entry is a linear combination of the $u_j$'s.  
$$\alpha \alpha'=\pmatrix{a_{11} & a_{12} \cr 
	     a_{21} & a_{22} \cr 
	     a_{31} & a_{32}} 
\pmatrix{0 & 1 \cr 
	 -1 & 0 } 
\pmatrix{a_{11} & a_{21} & a_{31} \cr 
	     a_{12} & a_{22} & a_{32} }$$ 
$$ = \pmatrix{ 0 & a_{11}a_{22}-a_{12}a_{21} & 
a_{11}a_{32}-a_{12}a_{31} \cr 
a_{21}a_{12}-a_{22}a_{11} & 0 & a_{21}a_{32}-a_{22}a_{31} \cr 
a_{31}a_{12}-a_{32}a_{11} & a_{31}a_{22}-a_{32}a_{21} & 0 } $$      
 is a $3\times 3$ matrix with each entry being a $7 \times 7$ block 
 matrix, which is a linear combination of the three linearly 
 independent matrices $B_j$ above. The condition that this composition 
 is zero says that the quadrics defined as the $2 \times 2$-minors of 
 the matrix $\alpha$, considered as quadrics in the $u_j$'s, have the 
 coefficients of $u_1^2$, $u_2^2$, $u_3^2$ respectively equal with the 
 coefficients of $u_0u_3$, $u_0u_1$, $u_0u_2$. 
\qed

(3.6) {\bf Proposition. } {\it If a matrix $\alpha$ with $\alpha \alpha'=0$ comes from 
an exact complex, then the three quadrics  given by its 
$2\times2$-minors are linearly independent.} 
 \medskip
{\it Proof.}Assume 
$$\alpha= \pmatrix{a_{11} & a_{12} \cr 
	     a_{21} & a_{22} \cr 
	     a_{31} & a_{32}}$$ 
is a matrix as above, with three linear dependent $2\times 2$-minors. 
 Without any loss of generality, may assume  
 
$$\alpha= \pmatrix{a_{11} & a_{12} \cr 
	     a_{21} & a_{22} \cr 
	     0 & 0} 
\hbox{ or \ \ } \alpha= \pmatrix{a_{11} & a_{12} \cr 
	                              0 & a_{22} \cr 
	                              0 & a_{32}}\ \ .$$ 
 
Indeed, after a linear change of the rows we may assume that 
the minor corresponding of the $2^{nd}$ and $3^{rd}$ row  is zero. 
A $2 \times 2$ determinant of linear forms is zero, iff either two rows or 
the 
two columns are linearly dependent. A further base change gives 
$\alpha$ 
the shape above. 
 
{\it Case 1.} If $\alpha$ has a zero row then the ideal contains a 
summand $V_4 \cO(-3)$, a contradiction. 
 
{\it Case 2.} In the second case we first note that  $a_{11} = l_0u_0+l_1u_1+l_2u_2+l_3u_3$ is 
non-zero. 
Indeed, otherwise  
$$ 0 \leftarrow I_A \leftarrow 3V_4 \cO(-3) \leftarrow 
S\Omega^3 
\leftarrow 0$$ 
would be exact, and $A$ could have codimension 2 at most. 
 
 Consider now  
$$V_4 \cO(-3) {\buildrel a_{11} \over \longleftarrow} S \Omega^3  
{\buildrel (u_0,u_1,u_2,u_3) \over \longleftarrow} 4V_1 \cO(-4).$$ 
The composition is the concatenation of the four $7 \times 7$ block matrices 
$l_1B_1+l_2B_2+l_3B_3 \mid l_0B_1-l_1B_3 \mid  
l_0B_2-l_2B_1 \mid l_0B_3-l_3B_2$.  
Precisely 2 blocks are linearly independent, because on one hand 
$a_{11}\cdot a_{22}=a_{11} \cdot a_{23} =0$, on the
 other hand the ideal $J$ in 
(4.1) has too few syzygies to allow an subideal ideal of type  
$a_{11} \cdot (a,b,c)$. So we get seven cubics with two skew symmetric $7 
\times 7$ matrices of relations. Since all 4 skew symmetric matrices 
have rang 6 in a general point of $\bP^6$ (eg. in the point  
$(1:2:3:4:5:6:7)$) unless they are identically zero, 
there is up to scalar a unique set of cubics 
whose relations they are: In each case the 7 principal pfaffians.  
These pfaffians are not proportional for two different blocks 
 for any values $(l_0,\ldots,l_3)$ unless the blocks themselves are 
proportional, as can be seen by a straight forward computation. 
 This is the desired contradiction. 
\qed

\bigskip 
 
{\bf 4. Moduli.} Denote by $X(1,7)^v$ the open set of abelian surfaces 
with a very ample polarization of class $(1,7)$. For each 
$A \in X(1,7)^v$ we choose a $G_7$-equivariant embedding 
$A \hookrightarrow \bP^6$. Its syzygy determine a $3 \times 2$ matrix 
$\alpha=\alpha_A$ as in Theorem 3.2. $\alpha$ is 
determined by $A$ up to conjugation with $GL(3,\bC) \times Sl(2,\bC)$. 
The ideal $I=I_A \subset S=\bC[u_0,u_1,u_2,u_3]$ of minors of $\alpha_A$ is 
uniquely determined by $A$. We denote by $C_A \subset \bP(U)$ the zero 
loci of $I_A$. 
\bigskip 
 
(4.1) {\bf Proposition.} {\it $C_A \subset \bP(U)$ is a projectively 
Cohen-Macaulay curve of degree $3$ and arithmetic genus $0$. } 
\medskip 
{\it Proof. } Denote by $S$ the graded ring $\bC [u_0,u_1,u_2,u_3] = 
S(U')= \oplus _{\ell \ge 0}S^\ell U'$. By the Hilbert-Burch Theorem (cf. [E],  
Thm 20.15) the complex 
$$ 0 \longleftarrow S/I  \longleftarrow S \longleftarrow \oplus_1^3 
S(-2) {\buildrel \alpha \over \longleftarrow} \oplus_1^2 S(-3) \longleftarrow 
0$$ 
is exact unless the three quadric minors of $\alpha$ have a common factor. 
Since the quadrics are linearly independent by Proposition 3.6, the 
second possibility occurs only if $S/I$ has syzygies 
$$\matrix{  
1 & - & - & -  \cr 
- & 3 & 3 & 1 \cr}.$$ 
However 
$$J:=((\Delta_1,\Delta_2,\Delta_3)^\bot) = (u_1u_2,u_2u_3,u_3u_1, 
u_1^2+u_0u_3,u_3^2+u_0u_2,u_2^2+u_0u_1,u_0^2)$$ 
has syzygies 
$$\matrix{  
1 & - & - & - & - \cr 
- & 7 & 8 & - & - \cr 
- & - & 3 & 8 & 3 \cr}$$ 
and $I \subset J$ by Proposition 3.5. So the second possibility cannot 
occur and the Hilbert-Burch complex is exact. \qed

(4.2) {\bf Corollary.} {\it $A$ is uniquely determined by $C_A$.} 
\medskip 
{\it Proof.} $C_A$ determines the Hilbert-Burch matrix $\alpha$ 
up to conjugation, which in turn determines $\alpha'$, $\beta'$ 
hence $I_A$. \qed

(4.3) The Hilbert scheme $Hilb_{3t+1}(\bP^3)$ has two components of dimension 12  
and 15 (cf. [PS]): 
$$Hilb_{3t+1}(\bP^3)=H_1 \cup H_2$$ 
with  general points of $H_1,H_2$ and the intersection $H_1 \cap H_2$ 
corresponding respectively to a twisted cubic, a plane cubic union a point or 
a plane nodal cubic with an embedded point at the node. 
For all $C \in H_1$, $h^0(\bP^3,{\cal I}_C)=3$. The morphism 
$$f\colon \matrix{ H_1 & \longrightarrow & \bG(3,H^0(\bP^3,\cO(2)) 
 \cr	  C & \mapsto & H^0(\bP^3,{\cal I}_C(2)) } $$ 
is birational onto its image $H \subset \bG(3,10)$, regular precisely on 
$H_1 - H_1 \cap H_2$, cf. [EPS]. All varieties 
$H,H_1,H_2,H_1\cap H_2,f(H_1 \cap H_2)$ are smooth.  
 
Consider 
$$H(\Delta):=H \cap \bG(3,J_2) \subset  \bG(3,H^0(\bP^3,\cO(2)).$$ 
Since $\bG(3,J_2)$ does not intersect $f(H_1 \cap H_2)$ we can regard 
$H(\Delta)$ as a subvariety of $H_1$ as well.  
$H(\Delta)$ has dimension at least 3 in every point by dimension count. 
 
 We are grateful to 
Geir Ellingsrud for pointing out to us, that such varieties were 
studied by Mukai. 
\medskip 
 
(4.4) {\bf Theorem.} {\it $H(\Delta)$ is a smooth prime Fano 3-fold of 
genus 12.} 
\medskip 
 
{\it Proof.} Mukai [Muk89,92] proves that $H(\delta) = H \cap 
\bG(3,\delta^\bot) 
\subset \bG(3,H^0(\bP^3,\cO(2))$ is a smooth prime Fano 3-fold 
for a general   net of quadrics $(\delta_1,\delta_2,\delta_3)$ in $\check 
\bP(U)$.  The proof that $\Delta_1,\Delta_2,\Delta_3$ is general in 
this sense, i.e. that $H(\Delta)$ is a smooth connected Fano 3-fold, is 
postponed until we have considered different models of $H(\Delta)$. 
 
\medskip 
 
(4.5) Consider on $L=\bC^7$ the net $\eta_{klein}\ \colon \Lambda^2L \to 
W'=\bC^3$ of alternating forms defined by the matrix 
$$\eta_{klein}\ = \pmatrix{ 0 & 0 & 0 & 0 & 0 & -y_1 & y_0 \cr 
0 & 0 & 0 & 0 & -y_2 & 0 & y_1 \cr 
0 & 0 & 0 & -y_0 & 0 & 0 & y_2 \cr 
0 & 0 & y_0 & 0 & y_1 & -y_2 & 0 \cr 
0 & y_2 & 0 & -y_1 & 0 & y_0 & 0 \cr 
y_1 & 0 & 0 & y_2 & -y_0 & 0 & 0 \cr 
-y_0 & -y_1 & -y_2 & 0 & 0 & 0 & 0 \cr 
}$$ 
and 
$$\bG(3,L,\eta_{klein})\ = \{ E \in \bG(3,L) \mid \Lambda^2 E \subset  
Ker (\eta_{klein}\colon \Lambda^2 L \to W' \}.$$ 
As zero loci of a section of a homogenous bundle on the Grassmanian, 
$\bG(3,L,\eta_{klein})$ is a prime Fano 3-fold of genus 12, if it is 
smooth of expected dimension. Smoothness follows from the criterion 
in section 1 of [Muk89] by computation. 
\medskip  
 
(4.6) Let $F = \{f=0\} \subset \bP^2$ be a plane quartic. The variety 
$$VSP(F,6)= \{ \{l_1,\ldots,l_6\} \in Hilb_6(\check \bP^2)  
\mid f = l_1^4+\ldots +l_6^4 \}$$ 
of sums powers presenting f was studied by Rosanes [Ros] 1873, Scorza
[Sco1,2]
and more recently  by Mukai [Muk89,92]. It is a prime Fano 
3-fold of genus 12 for general F.  Consider $f_{klein}\ = 
v_1^3v_2+v_2^3v_3+v_3^3v_1$ the well-known equation of the modular 
curve $\overline X(7) \subset \bP^2=\bP(W)$ due to Felix Klein, [K] \S 4. 
\medskip 
 
(4.7) {\bf Theorem.} {\it 
$$ H(\Delta) \cong \bG(3,L,\eta_{klein}\ ) \cong VSP(\overline X(7),6).$$ } 
 
{\it Proof.} Every prime Fano 3-fold $V_{22}$ of genus 12 (hence degree 22) 
has these 3 descriptions [Muk89,92] over an algebraically closed 
field 
of characteristic 0. That these special ones  
correspond to each other follows from [Schr], where the relation 
between the defining data and the isomorphism between the different
models is explained: 
 
$\eta_{klein}$ can be identified with the Tor-multiplication 
$$\Lambda^2 Tor_1^S(S/J,\bC)_2 \longrightarrow Tor^S_2(S/J,\bC)_4.$$ 
Note that these Tor-groups are 7 respectively 3-dimensional, 
cf. (4.1.1), and in fact $Tor_1^S(S/J,\bC)_2 =L$ and $Tor_2^S(S/J,\bC)_4=W'$,  
because the minimal resolution of $S/J=\bC\oplus U'\oplus W'$ over $S$ has  
the form: 
$$ 
0 \leftarrow S/J \leftarrow S \leftarrow LS(-2) \leftarrow M_1S(-3)\oplus  
W'S(-4)\leftarrow M_1S(-5)\leftarrow WS(-6)\leftarrow 0 
 \ \ . 
$$ 
On the other hand the ideal $I_{pfaff}$ generated by the 
$6 \times 6$ -Pfaffians of $\eta_{klein}$ gives a Gorenstein ring 
$A=\bC[y_0,y_1,y_2]/I_{pfaff}\ $ of codimension 3 cf.  
[BE]. A is artinian and the dual socle generator is $f_{klein}$. 
This completes the proof of Theorem 4.7 and 4.4. 
\qed 
 
(4.8) {\bf Remarks.} (1) The discriminant of the net of quadrics $\delta$ 
     is another quartic, which comes with a natural vanishing theta 
     characteristic, cf. 
 Scorza [1889,1899] and [DK]. In our case this is again the Klein 
     quartic. 
For general $\delta$ this is a different quartic than the dual socle 
     quartic, 
see [Schr] for more details. The fact that  
$$\{\hbox{\it quartics} \}\  - \to \{\hbox{\it quartics with an odd theta
characteristic} \}  $$
is birational over $\bC$ was discovered by Scorza. A more recent 
treatment is given in [DK], and with different point of view in [Schr].  
 
(2) If we take Mukai's results for granted, then it is clear that 
the quartic for the sum of powers has to coincide with the Klein 
quartic, because this curve is uniquely determined by its symmetry group: 
combine [ACGH] Ex. I F-17 with [H] Ex. IV 5.7 (b), or combine
[H] Ex. IV 5.7 (a) and the appendix.  
 
\medskip 
 
(4.9) {\bf Theorem.} {\it The Moduli space $X(1,7)$ is birational to 
$VSP(\overline X(7),6)$. } 
\medskip 
{\it Proof.} By 4.2 we have an immersion 
$$X(1,7)^v \hookrightarrow H(\Delta) \cong VSP(\overline X(7),6).$$ 
Since both varieties are irreducible and 3-dimensional the result
follows. \qed

(4.10) {\bf Theorem.} {\it $X(1,7)$ is rational with the rational 
map to $\bP^3$ defined over $\bQ$. } 
 
\medskip 
 
{\it Proof.} It suffices to prove that $\bG(3,L,\eta_{klein})$ is  
rational over $\bQ$.   
 
\noindent For a general point $p \in \bG(3,L,\eta)$ the triple projection 
defined by $| H -3p|$ defines a birational map 
$$\bG(3,L,\eta) - \to \bP^3,$$ 
(oral communication of Mukai). Its base loci consists of the 6 conics passing  
through $p$. 
This map is defined over $\bQ$, if the point is defined over $\bQ$. 
 
However the only readily visible rational point of
$\bG(3,V,\eta_{klein})$ 
is the point $p_e$ corresponding to the curve $C_e \subset \bP^3$  
defined by $u_1u_2=u_2u_3=u_3u_1=0$. Interpreted as a sum of powers this  
corresponds to the degenerate presentation 
$$ (v_1+\epsilon v_2)^4-v_1^4 + (v_2+\epsilon v_3)^4-v_2^4 +  
(v_3+\epsilon v_1)^4-v_3^4 = 4 \epsilon f_{klein}\ ,$$ 
viewed over $\bQ[\epsilon]/(\epsilon^2)$. For this reason we call $p_e$ 
the equational point. 
 
\medskip 
 
From the explicit form of $\eta=\eta_{klein}$ we see three lines  
$L_1,L_2,L_3 \subset \bG(3,L,\eta)$ passing through 
$p_e$. So $|H-3p_e|$ has larger dimension than for a general point. We pass to 
a subsystem. $|H-2L_1-2L_2-2L_3|$ has dimension 3. Its base loci consist 
of $L_1 \cup L_2 \cup L_3$ with each line with a 4-fold structure: 
The normal bundle of each line is $\cO_L \oplus \cO_L(-1)$. 
Hence ${\cal I}_L^2/{\cal I}_L^3$ has a summand  $\cO_L$ and this lies in 
the base loci since $|H-2L_1-2L_2-2L_3| \subset |H-3p_e|$. 
Note that the 4-fold structure is not generically a complete intersection.  
 
Two general hyperplanes $H_1,H_2 \in |H-2L_1-2L_2-2L_3|$ intersect along 
each line in a 5-fold structure. Hence the residual curve has degree 7. 
It intersects each line in one point. Hence a third general hyperplane 
$H_3 \in |H-2L_1-2L_2-2L_3|$ intersects the residual curve in these 3 points 
with multiplicity 2 and a single further point, i.e.  $|H-2L_1-2L_2-2L_3|$ 
defines a birational map to $\bP^3$.   
\medskip 
 
Instead of giving all the details of the arguments above we prefer to describe 
the inverse map 
$$ \psi \colon \bP^3 - \to \bG(3,V,\eta) = \bG(3,V) \cap \bP^{13}  
\subset \bP^{34} $$ 
explicitely.  
 
Consider the $3 \times 7$ matrix 
 
$$\Psi = \pmatrix{ -t_0t_3 & t_0t_1+t_2^2 & -t_3^2 & 0 & t_1t_3 & -t_2t_3 & 0 \cr 
t_1^2+t_0t_3 & -t_2^2 & -t_0t_2 & -t_1t_2 & 0 & t_2t_3 & 0 \cr 
t_0t_1^2+t_1t_2^2+t_0^2t_3 & t_2t_3^2 & t_1^2t_3+t_0t_3^2 & 0 & 0 & t_0t_2t_3 &   
t_1t_2t_3 \cr} $$ 
 
The rational map from $\bP^3$ with coordinates $t_0,\ldots,t_3$ to the Grassmanian 
$\bG(3,7)$ defined by $\Psi$ is for the given basis of $V$ the desired rational  
parametrization $\psi$. $\psi$ is bi-regular on $\bP^3-\{t_1t_2t_3=0\}$.  
\qed 
 
(4.11) {\bf Corollary.} {\it The rational universal family of $3 \times 2$ matrices 
on $\bP^3 \times \bP^3$ with coordinates $t,u$ is  
given by 
$$\alpha(t)=\pmatrix{ t_0u_1+t_2u_2 & -t_2u_0 & -t_1u_1 \cr 
	   	t_2u_2 & -t_0u_2-t_3u_3 & t_3u_0 \cr 
		u_1 & u_2 & u_3 \cr 
		t_1 & t_2 & t_3 \cr}.$$ }	 
 
\medskip 
 
{\it Proof.} The 3x3 minors of $\alpha(t)$ give 2 forms of bidegree (2,2) and two  
further forms of bidegree (3,2), (2,3) respectively. The 3 forms of degree 2 in u's 
are simply the product of the matrix of basis elements of $J$ with $^t\Psi$. 
The form of degree 3 in the u's is dependent unless $t_1=t_2=t_3=0$. Thus for any  
given point $[t_0:t_1:t_2:t_3] \ne [1:0:0:0]$ the minimal version of the matrix  
$\alpha(t)$ is a 3 x 2 matrix of linear forms in the u's, whose minors
are 
the  
desired 3 quadrics. 
\qed

Consider the vector 
$$D = (x_0x_3x_4,x_0x_1x_6,x_0x_2x_5,x_2^2x_3+x_4x_5^5,x_1^2x_5+x_2x_6^2,x_4^2x_6+ 
x_1x_3^2,x_1x_2x_4+x_3x_5x_6-x_0^3)$$ 
of $\tau$-invariant cubics in $S^3V_3$. 
 
\medskip 
 
(4.12) {\bf Corollary.} {\it The rational family on $\bP^3 \times \bP^6$ defined by  
the $H_7$-invariant 
subspace of cubics generated by the $\tau$-invariant forms 
$$(g_1 ,g_2, g_3)= D \cdot ^t\Psi$$ 
has as its fibres a dense subfamily of the universal family of $G_7$-invariant 
abelian surfaces of type (1,7). } 
\medskip 
{\it Proof.} Observe that, with the notations from A5. one has: 
$$J=(f_3,\ f_1,\ f_2,\ f_4,\ f_5,\ f_6,\ f_0)$$ 
and 
$$D=(f_3e_0,f_1e_0,f_2e_0,f_4e_0,f_5e_0,f_6e_0,f_0e_0)$$ 
The entries of $J$ and $D$ correspond to each other as 
 elements in  
two isomorphic irreducible $SL_2(\bZ_7)$-modules. 
\qed

{\bf Appendix} 
  
\medskip 
 
A1. {\bf The Heisenberg Group $H_7=H_{(1,7)}$ and the extended variant  
$G_7=H_7\rtimes \bZ_2$}    
We use notations similar to those from [HM] and [Man86], [Man89]

(A1.1) {\bf The direct construction of $H_{(1,7)}$ as a subgroup of $SL_7(\bC)$}

Let $V=Map(\bZ _7,\bC)$. On $V$ consider the automorphisms $\sigma $, $\tau $  
defined through: 
$$\eqalign{\sigma x(j)&=x(j+1) \cr 
\tau x(j)&=\varepsilon ^j x(j) \cr} 
$$ 
where $\varepsilon = {e ^ {2\pi i}\over {7}} \in \mu _7$ and let $H_7$ be the  
group generated by $\sigma$, $\tau$ (the Heisenberg group $H_{(1,7)}$). Then  
$H_7$ is generated by matrices $A_{ij}$  
of the form: 
$$ 
A_{j\ell} = \left( \varepsilon ^ {aj+b} \delta _{j,\ell +c} \right ) \ , 
\hbox{where}\  a,b,c \in \bZ _7 
$$ 
and  has order $343$. 
 
The Galois group $\Theta $ of $\bQ (\varepsilon )$ over $\bQ$ acts on H and let  
$\theta $ be the generator given by $\theta (\varepsilon ) = \varepsilon ^3$. Then  
$\theta ^3 = \hbox {\sl complex conjugation}$. The group $H$ is a central extension  
preserved by the action of $\Theta $: 
$$ 
1 \rightarrow \mu _ 7 \rightarrow H \rightarrow \bZ_7 \times \bZ _7 \rightarrow 0 
$$ 
(we identified tacitly $\mu _7$ and $\bZ_7$ in the last nonzero term and $\sigma  
\mapsto (1,0)$, 
$\tau \mapsto (0,1)$ ).

The irreducible $H$-module $V$ produces $5$ more by the composition with the  
automorphisms  
$\theta ^i \in \Theta $; denote by $V_i$ the representation 
$ 
H  {\buildrel \theta ^i \over \longrightarrow} H \rightarrow AutV$.  
These $6$ representations are inequivalent, as one sees computing their characters,  
and together with the characters of  
$\bZ_7 \times \bZ _7$ are all irreducible characters of $H$. 
 
Let $\Phi : \bZ_7 \times \bZ_7 \rightarrow H $ be the map  
$$\Phi (m,n) = \varepsilon ^{4mn} \sigma ^m \tau ^n $$ 
and consider  
$$ \omega :\mu _7 \times (\bZ_7 \times \bZ_7) \rightarrow  H \hbox{\quad \sl given by  
\quad}  
\omega (\alpha ,z) =\alpha \Phi (z) $$ 
Then $\omega $ is a bijection and the product in $H$ corresponds to  
$$ (\alpha, z)\cdot (\alpha ',z')=(\alpha \alpha ' B(z,z'),z+z') \hbox{\quad \sl where 
\quad } B(m,n;m',n')=\varepsilon ^{3(mn'-m'n)}$$ 
Let $N$ be the normaliser of $H$ in $SL_7(\bC)$. Then $N$ is a central extension: 
$$ \displaylines{ 
1\rightarrow H \rightarrow N {\buildrel \alpha \over \longrightarrow} SL_2(\bZ_7)  
\rightarrow 1\cr 
\hbox{\sl where } \alpha (x) = \hbox{\sl the automorphism of $\bZ_7 \times \bZ_7$  
preserving B, induced by conjugation by } x \in Aut(H).} 
$$ 
Note that for $u \in SL_2(\bZ_7)$ the action $\gamma _u$ on $H$ can be expressed: 
$$ 
\gamma _u \omega (\alpha , z) = \omega (\alpha , u(z)). 
$$ 
We have seen that any polarized abelian surface $(A,L)$ with a polarization of type  
$(1,7)$, $L$ very ample and a fixed 
level structure is canonically embedded in $\bP ^6 :=\hbox {\sl the projective space  
of lines in } V$. Moreover, we may suppose that $L$ is symmetric with respect to the  
origin 
of $A$, because, by a change of the origin, we can realize this situation. Then the  
map $x \mapsto -x$ on $A$ extends to an automorphism of order $2$ of $\bP^6$,  
induced by  
$\iota \in SL(V)$, $\iota x(j)=-x(-j)$. Therefore we consider from now on abelian  
surfaces in $\bP^6$ invariant under the action of $G=H \rtimes \bZ_2$. 
\medskip
{\bf Remark}. If $\{e_j\}_j$ is the canonical basis of $V=Map(\bZ_7,\bC)$, i.e.  
$e_j(\ell)=\delta _{j\ell}$, and $\{x_j\}_j$ is the dual basis of $V\dual=V_3$, then  
the action of $\sigma , \tau $ on $V$ and on $H^0(\bP^6, \cO(1))=V\dual =V_3$ is given  
by:

$$ 
\vbox{ 
\halign { $#$ \hfil  &  $ # $\hfil   & \quad \quad $#$ \hfil &   $# $\hfil  \cr 
& & & \cr 
\sigma e_j  &  = e_{j-1} &  \sigma x_j  & = x_{j-1}  \cr 
\tau e_ j   &  =\varepsilon ^j e_j   &  \tau x_j  & = \varepsilon ^{-j} x_j\cr} 
} 
$$ 
\medskip 
 
(A1.2) {\bf Character table of $G$} 
 (cf. [Man86] for $H_5\rtimes \bZ_2$): 
 
\medskip 
 
\vbox 
{\tabskip=-1.7pt \offinterlineskip 
\hrule 
\halign to 300pt 
{ 
\strut # 
& \vrule  # \tabskip=1em plus2em  
& \hfil # \hfil 
& \vrule #  
& \hfil # \hfil 
& \vrule # 
& \hfil # \hfil  
& \vrule # 
& \hfil # \hfil 
& \vrule # \tabskip=0pt  
\cr 
\noalign{\hrule} 
&& $\{\alpha \}$  
&& $C_{m,n}$  
&& $C_\alpha $ 
&&  
&\cr 
\noalign{\hrule} 
&& $1$ 
&& $1$ 
&& $1$ 
&& $I$ 
&\cr 
\noalign{\hrule} 
&& $7\theta^i(\alpha )$ 
&& $0$ 
&& $\theta^i (\alpha ) $ 
&& $V_i$ 
&\cr 
\noalign{\hrule} 
&& $1$ 
&& $1$ 
&& $-1$ 
&& $S$ 
&\cr 
\noalign{\hrule} 
&& $7\theta^i(\alpha )$ 
&& $0$ 
&& $-\theta^i (\alpha ) $ 
&& $V_i^\sharp$ 
&\cr 
\noalign{\hrule} 
&& $2$ 
&& $\varepsilon ^{sm+tn}+\varepsilon ^{-sm-tn}$ 
&& $0$ 
&& $Z_{s,t}$ 
&\cr 
\noalign{\hrule} 
\cr} 
} 
 
\medskip 
where:\smallskip 
 
$\{\alpha \}$ is the conjugacy class containing only the central element 
$\alpha \in \mu _7$, \smallskip  
 
$C_{m,n} = \{ (\alpha ,m,n),(\alpha ,-m,-n) | \alpha \in \mu _7\} 
(\hbox {\sl \quad with $(m,n)\not= (0,0)$})$ \smallskip 
 
and $C_\alpha =\{ (\alpha,m,n)\iota \ 
| m,n \in \bZ_7 \}$.\smallskip  
 
Thus there are $7$ classes $\{\alpha \}$, $24$ classes $C_{m.n}$ (each with 14  
elements) and $7$ classes $C_\alpha $ (each with $49$ elements). We denote by $Z$  
the sum of all 
$24$ $Z_{s,t}$. 
 
\medskip 
 
(A1.3) {\bf Useful formulae}

We have the following formulae: 
$$ 
\eqalign{ 
V_i \otimes V_i     & = 3V_{i+2} \oplus 4V_{i+2}^\sharp \cr 
V_i \otimes V_{i+1} & = 3V_{i+4} \oplus 4V_{i+4}^\sharp \cr 
V_i \otimes V_{i+2} & = 3V_{i+1} \oplus 4V_{i+1}^\sharp \cr 
V_i \otimes V_{i+3} & = I\oplus Z \cr 
}\quad 
\eqalign{ 
\wedge ^2 V_i     & = 3V_{i+2} \cr 
\wedge ^3 V_i     & = V_{i+1}\oplus 4V_{i+1}^\sharp \cr 
\wedge ^4 V_i     & = V_{i+4}\oplus 4V_{i+4}^\sharp \cr 
\wedge ^5 V_i     & = 3V_{i+5} \cr 
}\quad 
\eqalign{ 
\wedge ^6 V_i     & = V_{i+3} \cr 
\wedge ^7 V_i     & = I \cr 
\cr 
\cr 
} 
$$

$$ 
\eqalign{ 
S ^2 V_i     & = 4V_{i+2}^\sharp \cr 
S ^3 V_i     & = 8V_{i+1}\oplus 4V_{i+1}^\sharp \cr 
S ^4 V_i     & = 10V_{i+4}\oplus 20V_{i+4}^\sharp \cr 
S ^5 V_i     & = 38V_{i+5}\oplus 28V_{i+5}^\sharp \cr 
S ^6 V_i     & = 56V_{i+3} \oplus 76V_{i+3}^\sharp \cr 
S ^7 V_i     & = 8I \oplus 28S\oplus 35Z \cr 
S ^8 V_i     & = 197V_i \oplus 232 V_i ^\sharp \cr 
S ^9 V_i     & = 375V_{i+2} \oplus 340 V_{i+2} ^\sharp \cr 
}\quad 
\eqalign{ 
S ^{10} V_i     & = 544V_{i+1} \oplus 600 V_{i+1} ^\sharp \cr 
S ^{11} V_i     & = 908V_{i+4} \oplus 852 V_{i+4} ^\sharp \cr 
S ^{12} V_i     & = {1\over 7}\left( {1\over 2}{18\choose 6}-42\right)V_{i+5} \oplus 
{1\over 7}\left( {1\over 2}{18\choose 6}+42\right)V_{i+5}^\sharp \cr 
S ^{13} V_i     & = {1\over 7}\left( {1\over 2}{19\choose 6}-42\right)V_{i+3} \oplus 
{1\over 7}\left( {1\over 2}{19\choose 6}+42\right)V_{i+3}^\sharp \cr 
S ^{14} V_i     & = 12\cdot384 I\oplus 12\cdot 374S \oplus 48\cdot 618 Z \cr 
etc. & \cr 
\cr 
} 
$$

$$ 
\eqalign{ 
H^0(\Omega ^3(3)) & = 0 \cr 
H^0(\Omega ^3(4)) & = \wedge ^3V=V_1\oplus 4V_1^\sharp\cr 
H^0(\Omega ^3(5)) & =16V_2\oplus 16V_2^\sharp\cr 
H^0(\Omega ^3(6)) & = 56V\oplus 64V^\sharp\cr 
H^0(\Omega ^3(7)) & = 24I \oplus 24S \oplus 49Z \cr 
} 
\quad 
\eqalign{ 
H^0(\Omega ^3(8)) & = 405V_3\oplus 420V_3^\sharp\cr 
H^0(\Omega ^3(9)) & = 880V_5\oplus 880V_5^\sharp\cr 
H^0(\Omega ^3(10)) & =1704V_4\oplus 1728V_4^\sharp\cr 
etc. & \cr 
} 
$$

Observing that the trace of $\iota $ on $H^0({\cal O }_A(k))$ is $-1$ for $k$ odd  
and $4$ for $k$ even, one deduces:

$$ 
\eqalign{ 
H^0({\cal O }_A (1)) & = V_3 \cr 
H^0({\cal O }_A (2)) & = 4V_5 ^\sharp \cr 
H^0({\cal O }_A (3)) & = 5V_4\oplus 4V_4 ^\sharp \cr 
H^0({\cal O }_A (4)) & = 6V_1\oplus  10 V_1^\sharp \cr 
H^0({\cal O }_A (5)) & = 13V_2 \oplus 12V_2^\sharp \cr 
H^0({\cal O }_A (6)) & = 16V \oplus 20V ^\sharp \cr 
H^0({\cal O }_A (7)) & = 3I \oplus 4S \oplus 7Z \cr 
\cr 
} 
\quad 
\eqalign{ 
H^0({\cal O }_A (8)) & = 30V_3 \oplus 34V_3^\sharp \cr 
H^0({\cal O }_A (9)) & = 41V_5\oplus 40V_5^\sharp \cr 
H^0({\cal O }_A (10)) & = 48V_4 \oplus 52V_4^\sharp \cr 
H^0({\cal O }_A (11)) & = 61V_1 \oplus 60V_1^\sharp \cr 
H^0({\cal O }_A (12)) & = 70V_2 \oplus 74V_2^\sharp \cr 
H^0({\cal O }_A (13)) & = 85V \oplus 84V^\sharp \cr 
H^0({\cal O }_A (14)) & = 16I\oplus 12S \oplus 28Z \cr 
etc. & \cr 
} 
$$ 
 
\medskip 
 
A2 {\bf The Normaliser $N$ of $H$ in $SL_7(\bC)$.} 
Via the map $\alpha : N \to SL_2(\bZ_7)$ we get, entirely like in [HM], 
that $N$ is a semidirect product $H\rtimes SL_2(\bZ_7)$. Then, in fact $N \subset  
SL_7(\bQ(\varepsilon))$. Thus $\Theta $ acts on $N$ and all $V_i$ are also  
$N$-modules.  
One shows, like in [HM], that $V_i \otimes V_i\dual \cong I\oplus Z$, for  
all  
$i$, where $Z$ is the space of trace $0$, and as a $\bZ_7 \times \bZ_7$-module is  
the sum of all $48$ nontrivial modules. In fact, as a $N/\mu_7$-module it is  
irreducible and its character takes values in $\bQ$. 
  
\medskip 
 
A3 {\bf The Table of Characters of the Group $SL_2(\bZ_7)$ and their  
Multiplication.}  
First of all we make some notations: 
 
$$ 
\eqalign{ 
\alpha   = i\sqrt{7}\quad \quad \quad \quad \lambda _1  = \varepsilon ^{1^2}- 
\varepsilon ^{-1^2} = \varepsilon -\varepsilon ^6  
 & \quad \quad \quad \eta _1 = \varepsilon +\varepsilon ^6 \cr 
\alpha ^+  = {1 \over 2}(1+\alpha ) \quad \quad \lambda _2  =  
\varepsilon ^{2^2}-\varepsilon ^{-2^2} =\varepsilon ^4-\varepsilon ^3 &  
\quad \quad \quad \eta  _2 = \varepsilon ^4 +\varepsilon ^3 \cr 
\alpha ^-  = {1 \over 2}(1-\alpha ) \quad \quad \lambda _3  =  
\varepsilon ^{3^2}-\varepsilon ^{-3^3} =\varepsilon ^2-\varepsilon ^5 & \quad \quad  
\quad  \eta _3 = \varepsilon ^2 +\varepsilon ^5 \cr 
} 
$$ 

Then we have the following equalities, useful in computations: 
 
$$ 
\eqalign{ 
\varepsilon +\varepsilon ^2 +\varepsilon ^4  = -\alpha ^- \cr 
\varepsilon ^3+\varepsilon ^5 +\varepsilon ^6   = -\alpha ^+  
}\quad \quad 
\eqalign{ 
\lambda _1 +\lambda _2 +\lambda _3 = \alpha \cr 
 \lambda _1\lambda _2\lambda _3 = \alpha \cr 
}\quad \quad 
\eqalign{ 
\eta _1 +\eta _2 +\eta _3 =-1 \cr 
 \eta _1\eta _2\eta _3=1 \cr 
} 
$$ 
and 
$$  
\eqalign{ 
\lambda _1 ^2=\eta _3 -2 \cr 
\lambda _2 ^2=\eta _1 -2 \cr 
\lambda _3 ^2=\eta _2 -2 \cr 
}\quad \quad 
\eqalign{ 
\lambda _1 \lambda _2 =\eta _3 -\eta _2 \cr 
\lambda _2 \lambda _3 =\eta _1 -\eta _3 \cr 
\lambda _3 \lambda _1 =\eta _2 -\eta _1 \cr 
}\quad \quad 
\eqalign{ 
\alpha \eta _1 = \lambda _1-2\lambda _2\cr 
\alpha \eta _2 = \lambda _2-2\lambda _3\cr 
\alpha \eta _3 = \lambda _3-2\lambda _1\cr 
} 
$$ 
 
The general shape of the elements in $N$ is:
\medskip 
\vbox{ 
\halign {\hfil # &  # \hfil \cr 
$A_{jk}$ &$= \pm {1 \over \sqrt 7}\varepsilon ^{aj^2+bjk+ck^2+dj+ek+f} \quad \quad  
(a,b, \ldots ,f \in \bZ_7, \ b\not = 0 )$ \cr 
& \cr 
$A_{jk}$ &$= \pm \varepsilon ^{aj^2+bj+c} \delta_{j,dk+e}\quad \quad (a,b, \ldots ,  
e \in \bZ_7, \ d\not = 0 )$\cr 
} 
} 
\medskip 
 
\noindent where the signs are chosen to have $det (A_{jk}) = 1$. 
 
For the convenience of the computations, it is useful to identify some elements in  
$N=H\cdot 
SL_2(\bZ_7)$ and their images in $SL_2(\bZ_7)$:

$$ 
\eqalign{ 
\mu x(j) &=x(2j) \quad ({\rm resp.\ } \mu e_j = e_{j/2}) \quad \hbox {\sl  
corresponding \ to \ } \overline {\mu} = \pmatrix{2 & 0 \cr 
                                                             0 & 4 } \in SL_2(\bZ_7)  
                                                             \cr 
\nu x(j) &=\varepsilon ^{j^2}x(j) \quad ({\rm resp.\ } \nu e_j = 
\varepsilon ^{j^2}e_j) \quad \hbox {\sl corresponding \ to \ } \overline {\nu} = 
\pmatrix{1 & 0 \cr 
         2 & 1    } \in SL_2(\bZ_7) \cr 
\delta x(j) &={i\over \sqrt{7}}\sum_k\varepsilon ^{kj}x(k) \quad ({\rm resp.\ }  
\delta e_j ={i\over \sqrt{7}}\sum_k\varepsilon ^{kj}e_k ) \quad {\sl corresponding  
\ to \ }  
\overline {\delta} = 
\pmatrix{0 & -1 \cr 
         1 &  0    } \in SL_2(\bZ_7) \cr 
} 
$$ 
 
Observe that $\delta ^2 = \iota$ and that the elements in $SL_2(\bZ_7)$ are given  
according to: 
 
$$ 
\halign 
{\indent \quad # \hfil & \quad # \hfil & \quad # \hfil & \quad # \hfil \cr 
$\mu\sigma\mu^{-1}=\sigma^2$ & $\iota\sigma\iota = \sigma^{-1}$ & $\nu\sigma\nu^{-1} 
=\varepsilon ^{4\cdot 1\cdot 2}\sigma\tau^2$ & $\delta\sigma \delta^{-1}=\tau $ \cr 
$\mu\tau\mu^{-1}=\tau^4$ & $\iota\tau\iota = \tau^{-1}$ & $\nu\tau\nu ^{-1}=\tau$ &  
$\delta\tau\delta^{-1}= 
\sigma ^{-1}$\cr 
} 
$$

\vbox 
{\tabskip=-1.65pt \offinterlineskip 
\halign to 405pt 
{ 
\strut # 
& \vrule  # \tabskip=5pt plus 1pt 
& \hfil # \hfil 
& \vrule #  
& \hfil #  
& \vrule # 
& \hfil #   
& \vrule # 
& \hfil #  
& \vrule #  
& \hfil #  
& \vrule # 
& \hfil #   
& \vrule #  
& \hfil #  
& \vrule # 
& \hfil #  
& \vrule # 
& \hfil #  
& \vrule #  
& \hfil #  
& \vrule # 
& \hfil #  
& \vrule # 
& \hfil #  
& \vrule # \tabskip=0pt  
\cr 
\multispan{25} \hfil Character Table of $SL_2(\bZ_7)$ \hfil \cr  
\noalign{\bigskip \hrule} 
&&  
&& $id$ 
&& $\iota$ 
&& $\mu$ 
&& $\iota\mu$ 
&& $\nu$ 
&& $\nu ^3$ 
&& $\iota\nu ^3$ 
&& $\iota\nu$ 
&& $\delta$ 
&&  
&&   
&\cr 
&&  
&&  
&&  
&&  
&&  
&&  
&&  
&&  
&&  
&&  
&&  
&&   
&\cr 
&&  
&& $id$ 
&& $-id$ 
&& $\overline{\mu}$ 
&& $\overline{\iota\mu}$ 
&& $\overline{\nu}$ 
&& $\overline{\nu } ^3$ 
&& $\overline{\iota \nu }^3$ 
&& $\overline{\iota \nu }$ 
&& $\overline{\delta }$ 
&& ${22}\choose{52}$ 
&& ${52}\choose{55}$  
&\cr 
&&  
&&  
&&  
&&  
&&  
&&  
&&  
&&  
&&  
&&  
&&  
&&   
&\cr 
&&  
&& $1$ 
&& $1$ 
&& $56$ 
&& $56$ 
&& $24$ 
&& $24$ 
&& $24$ 
&& $24$ 
&& $42$ 
&& $42$ 
&& $42$  
&\cr 
\noalign{\hrule} 
&& $I$ 
&& $1$ 
&& $1$ 
&& $1$ 
&& $1$ 
&& $1$ 
&& $1$ 
&& $1$ 
&& $1$ 
&& $1$ 
&& $1$ 
&& $1$ 
&\cr 
\noalign{\hrule} 
&& $M_1$ 
&& $8$ 
&& $-8$ 
&& $-1$ 
&& $1$ 
&& $1$ 
&& $1$ 
&& $-1$ 
&& $-1$ 
&& $0$ 
&& $0$ 
&& $0$ 
&\cr 
\noalign{\hrule} 
&& $M_2$ 
&& $8$ 
&& $8$ 
&& $-1$ 
&& $-1$ 
&& $1$ 
&& $1$ 
&& $1$ 
&& $1$ 
&& $0$ 
&& $0$ 
&& $0$ 
&\cr 
\noalign{\hrule} 
&& $L$ 
&& $7$ 
&& $7$ 
&& $1$ 
&& $1$ 
&& $0$ 
&& $0$ 
&& $0$ 
&& $0$ 
&& $-1$ 
&& $-1$ 
&& $-1$ 
&\cr 
\noalign{\hrule} 
&& $U$ 
&& $4$ 
&& $-4$ 
&& $1$ 
&& $-1$ 
&& $\alpha $\rlap{$^-$} 
&& $\alpha $\rlap{$^+$} 
&& $-\alpha $\rlap{$^+$} 
&& $-\alpha $\rlap{$^-$} 
&& $0$ 
&& $0$ 
&& $0$ 
&\cr 
\noalign{\hrule} 
&&  $ \! U'=U\dual \! $  
&& $4$ 
&& $-4$ 
&& $1$ 
&& $-1$ 
&& $\alpha $\rlap{$^+$} 
&& $\alpha $\rlap{$^-$} 
&& $-\alpha $\rlap{$^-$} 
&& $-\alpha $\rlap{$^+$} 
&& $0$ 
&& $0$ 
&& $0$ 
&\cr 
\noalign{\hrule} 
&& $T_1$ 
&& $6$ 
&& $-6$ 
&& $0$ 
&& $0$ 
&& $-1$ 
&& $-1$ 
&& $1$ 
&& $1$ 
&& $0$ 
&& $\sqrt{2}$ 
&& $-\sqrt{2}$ 
&\cr 
\noalign{\hrule} 
&& $T_2$ 
&& $6$ 
&& $-6$ 
&& $0$ 
&& $0$ 
&& $-1$ 
&& $-1$ 
&& $1$ 
&& $1$ 
&& $0$ 
&& $-\sqrt{2}$ 
&& $\sqrt{2}$ 
&\cr 
\noalign{\hrule} 
&& $T$ 
&& $6$ 
&& $6$ 
&& $0$ 
&& $0$ 
&& $-1$ 
&& $-1$ 
&& $-1$ 
&& $-1$ 
&& $2$ 
&& $0$ 
&& $0$ 
&\cr 
\noalign{\hrule} 
&& $W$ 
&& $3$ 
&& $3$ 
&& $0$ 
&& $0$ 
&& $-\alpha $\rlap{$^+$} 
&& $-\alpha $\rlap{$^-$} 
&& $-\alpha $\rlap{$^-$} 
&& $-\alpha $\rlap{$^ +$} 
&& $-1$ 
&& $1$ 
&& $1$ 
&\cr 
\noalign{\hrule} 
&& $\! W'=W\dual \!$ 
&& $3$ 
&& $3$ 
&& $0$ 
&& $0$ 
&& $-\alpha $\rlap{$^-$} 
&& $-\alpha $\rlap{$^+$} 
&& $-\alpha $\rlap{$^+$} 
&& $-\alpha $\rlap{$^-$} 
&& $-1$ 
&& $1$ 
&& $1$ 
&\cr 
\noalign{\hrule} 
\cr} 
}

We indicate here also the multiplication table of the characters of
$SL_2(\bZ_7)$:
   
\vbox{ 
\halign{ $#$ \hfil &  $#$ \hfil & \quad  \quad $#$ \hfil & $#$ \hfil \cr  
M_1\otimes M_1 &  = I \oplus 3M_2 \oplus 3L \oplus 2T \oplus W \oplus W'   &  \cr 
M_1\otimes M_2 &  = 3M_1 \oplus 2U \oplus  2U'\oplus 2T_1 \oplus 2T_2 &  
M_2\otimes M_2 &  = I \oplus 3M_2 \oplus 3L \oplus 2T \oplus W \oplus W'  \cr 
M_1\otimes L   &  = 3M_1 \oplus U \oplus U' \oplus 2T_1 \oplus 2T_2 & 
M_2\otimes L   &  = 3M_2 \oplus 2L \oplus 2T \oplus W \oplus W' \cr 
M_1\otimes U   &  = 2M_2 \oplus L \oplus T \oplus W & 
M_2\otimes U   &  = 2M_1 \oplus U \oplus T_1 \oplus T_2   \cr 
M_1\otimes U'  &  = 2M_2 \oplus L \oplus T \oplus W'& 
M_2\otimes U'  &  = 2M_1 \oplus U' \oplus T_1 \oplus T_2  \cr 
M_1\otimes T_1 &  = 2M_2 \oplus 2L \oplus 2T \oplus W \oplus W' & 
M_2\otimes T_1 &  = 2M_1 \oplus U \oplus U' \oplus 2T_1 \oplus 2T_2  \cr 
M_1\otimes T_2 &  = 2M_2 \oplus 2L \oplus 2T \oplus W \oplus W' & 
M_2\otimes T_2 &  = 2M_1 \oplus U \oplus U' \oplus 2T_1 \oplus 2T_2   \cr 
M_1\otimes T   &  = 2M_1 \oplus U \oplus U' \oplus 2T_1 \oplus 2T_2  & 
M_2\otimes T   &  = 2M_2 \oplus 2L \oplus 2T \oplus W \oplus W'   \cr 
M_1\otimes W   &  = M_1 \oplus U \oplus T_1 \oplus T_2  & 
M_2\otimes W   &  = M_2 \oplus L \oplus T\oplus W    \cr 
M_1\otimes W'  &  = M_1 \oplus U'\oplus T_1 \oplus T_2 & 
M_2\otimes W'  &  = M_2 \oplus L \oplus T \oplus W'  \cr 
& & \cr} 
} 
\vbox{ 
\halign{ $#$ \hfil &  $#$ \hfil & \quad  \quad $#$ \hfil & $#$ \hfil \cr  
L\otimes L   &  = I\oplus 2M_2 \oplus 2L \oplus 2T \oplus W \oplus W'  &  \cr 
L\otimes U   &  = M_1 \oplus U \oplus U' \oplus T_1 \oplus T_2   &  \cr 
L\otimes U'  &  = M_1 \oplus U \oplus U' \oplus T_1 \oplus T_2   &  \cr 
L\otimes T_1 &  = 2M_1 \oplus U \oplus U' \oplus T_1 \oplus 2T_2   &  \cr 
L\otimes T_2 &  = 2M_1 \oplus U \oplus U' \oplus 2T_1 \oplus T_2   &  \cr 
L\otimes T   &  = 2M_2 \oplus 2L \oplus T \oplus W  \oplus W'  &  \cr 
L\otimes W   &  = M_2 \oplus L \oplus T  &  \cr 
L\otimes W'  &  = M_2 \oplus L \oplus T  &  \cr 
& &  \cr} 
} 
\vbox{ 
\halign{ $#$ \hfil &  $#$ \hfil & \quad  \quad $#$ \hfil & $#$ \hfil \cr  
U\otimes U    &  = L \oplus T \oplus W &  & \cr 
U\otimes U'   &  = I \oplus M_2 \oplus L &  
U'\otimes U'  &  = L \oplus T \oplus W'  \cr 
U\otimes T_1  &  = M_2 \oplus L \oplus T \oplus W' &  
U'\otimes T_1 &  = M_2 \oplus L \oplus T \oplus W  \cr 
U\otimes T_2  &  = M_2 \oplus L \oplus T \oplus W' &  
U'\otimes T_2 &  = M_2 \oplus L \oplus T \oplus W  \cr 
U\otimes T    &  = M_1 \oplus U' \oplus T_1 \oplus T_2 &  
U'\otimes T   &  = M_1 \oplus U \oplus T_1 \oplus T_2  \cr 
U\otimes W    &  = T_1 \oplus T_2 & 
U'\otimes W   &  = M_1 \oplus U \cr 
U\otimes W'   &  = M_1 \oplus U' & 
U'\otimes W'  &  = T_1 \oplus T_2 \cr 
& & \cr} 
} 
\vbox{ 
\halign{ $#$ \hfil &  $#$ \hfil & \quad  \quad $#$ \hfil & $#$ \hfil \cr  
T_1\otimes T_1 &  = I \oplus 2M_2 \oplus L \oplus T \oplus W \oplus W'   &  \cr 
T_1\otimes T_2 &  = 2M_2 \oplus 2L \oplus T  &  
T_2\otimes T_2 &  = I \oplus 2M_2 \oplus L \oplus T \oplus W \oplus W'  \cr 
T_1\otimes T   &  = 2M_1 \oplus U \oplus U' \oplus T_1 \oplus T_2  &  
T_2\otimes T   &  = 2M_1 \oplus U \oplus U' \oplus T_1 \oplus T_2  \cr 
T_1\otimes W   &  = M_1 \oplus U' \oplus T_1  &  
T_2\otimes W   &  = M_1 \oplus U' \oplus T_2    \cr 
T_1\otimes W'  &  = M_1 \oplus U \oplus T_1  &  
T_2\otimes W'  &  = M_1 \oplus U \oplus T_2    \cr 
& & \cr 
T\otimes T     &  = I \oplus 2M_2 \oplus L \oplus 2T   &   
W\otimes W     &  = T \oplus W'  \cr 
T\otimes W     &  = M_2 \oplus L \oplus W' & 
W\otimes W'    &  = I \oplus M_2  \cr 
T\otimes W'    &  = M_2 \oplus L \oplus W & 
W'\otimes W'   &  = T \oplus W  \cr} 
} 
 
\medskip  
 A4 {\bf Multiplications, exterior, symmetric powers of  
$N$ -- representations, some $H^0(\Omega ^3(j))$' s.} 
 
\medskip
 
\vbox{ 
\halign{\quad \quad \quad \quad $#$ \hfil   &  $#$ \hfil    & \quad  \quad $#$ \hfil  
& $#$ \hfil \cr V_{2j}\otimes V_{2j} &  = (U' \oplus W') \otimes V_{2j+2} & V_{2j} 
\otimes V_{2j+1} &=(U\oplus W ) 
\otimes V_{2j+4}  \cr 
V_{2j+1}\otimes V_{2j+1} &  = (U\oplus W) \otimes V_{2j+3} & V_{2j+1}\otimes V_{2j+2}  
&=(U'\oplus  
W' )\otimes V_{2j+5}  \cr 
& & &  \cr 
V_{2j}\otimes V_{2j+2} &  = (U\oplus W) \otimes V_{2j+1} & V_j\otimes V_{j+3} & 
=I\oplus Z   \cr 
V_{2j+1}\otimes V_{2j+3} &  = (U' \oplus W') \otimes V_{2j+2} & & \cr} 
}

\medskip  
 
\vbox{ 
\halign{\quad \quad \quad \quad $#$ \hfil   &  $#$ \hfil    & \quad  \quad $#$ \hfil  
& $#$ \hfil & \quad \quad $#$ & $#$ \hfil \cr  
\Lambda ^2 V_{2j} &  =  W' \otimes V_{2j+2} & \Lambda ^3 V_{2j} &=(I\oplus U') 
\otimes V_{2j+1} & \Lambda ^4 V_{2j} &= (I\oplus U)\otimes V_{2j+4}\cr 
\Lambda ^5 V_{2j} &  =  W \otimes V_{2j+5} &  &  &\cr 
& & & & &\cr 
\Lambda ^2 V_{2j+1} &  =  W \otimes V_{2j+3} & \Lambda ^3 V_{2j+1} &=(I\oplus U) 
\otimes V_{2j+2} & \Lambda ^4 V_{2j+1} &= (I\oplus U')\otimes V_{2j+5}\cr 
\Lambda ^5 V_{2j+1} &  =  W' \otimes V_{2j+6} &  &  &\cr 
& & & & & \cr  
\Lambda ^6 V_j &= V_{j+3} & \Lambda ^7 V_j &= I & & \cr 
& & & & &\cr 
} 
}

\medskip 
\vbox{ 
\halign{\quad \quad \quad \quad $#$ \hfil   &  $#$ \hfil \cr 
S ^2 V_{2j} &  =  U' \otimes V_{2j+2} \cr 
S ^3 V_{2j} &=(I\oplus L\oplus U)\otimes V_{2j+1} \cr 
S ^4 V_{2j} &= (L\oplus W\oplus U\oplus U'\oplus T_1\oplus T_2)\otimes V_{2j+4}\cr 
S^5 V_{2j} &  =  (I\oplus  M_1\oplus M_2\oplus 2L\oplus U\oplus U'\oplus T_1 
\oplus T_2\oplus 2T\oplus W')\otimes V_{2j+5} \cr 
etc. &\cr 
} 
} 
 
\medskip
\vbox{ 
\halign{\quad \quad \quad \quad $#$ \hfil   &  $#$ \hfil \cr 
H^0(\Omega ^3(3)) &  =  0 \cr 
H^0(\Omega ^3(4)) &  =  (I\oplus U')\otimes V_1 \cr 
H^0(\Omega ^3(5)) &  =  (L\oplus U' \oplus W' \oplus T_1 \oplus T_2 
\oplus T)\otimes V_2 \cr 
H^0(\Omega ^3(6)) &  =  (M_1\oplus M_2\oplus 3L \oplus 2U \oplus 2W \oplus W'  
\oplus 4T_1 \oplus 4T_2 \oplus 3T)\otimes V\cr 
H^0(\Omega ^3(7)) &  =  (I \oplus 2L \oplus U \oplus 2U' \oplus W' \oplus T_1 
\oplus T_2 \oplus T)(I\oplus Z) \oplus Z\cr 
etc. &\cr 
} 
}

A5 {\bf Concrete Decompositions of Certain $N$ or $SL_2(\bZ_7)$ 
Representations.}
  
Consider now the decomposition of $V$ into eigenspaces of $\iota $: 
$$ 
\halign 
{ \indent \quad # &  # \hfil  \cr 
$V$    & $ = V^+ \oplus V^-$ \quad where: \cr 
&  \cr 
$V^+$  & $ = span \ \{ e_{1^2}-e_{-1^2},\ e_{2^2}-e_{-2^2},\ e_{3^2}-e_{-3^2}\}  
        = span \ \{e_1-e_6, \ e_4-e_3, \  e_2-e_5\}$    \cr 
&  \cr 
$V^-$  & $ = span \ \{2e_0, \  e_1+e_6,\  e_4+e_3,\  e_2+e_5\} $ \cr 
} 
$$ 
Restricting $\mu $, $\nu $, $\delta $ to $V^+$ and $V^-$ respectively,  
one gets:
$$ 
\halign{\indent # & # & # & # & # & #\cr 
$\mu ^+$ & $= \pmatrix {0 & 0 & 1 \cr 
                    1 & 0 & 0 \cr 
                     0 & 1 & 0 },\quad $ & $\nu ^+$ & $= diag\ (\varepsilon ,  
\varepsilon ^2 , \varepsilon ^4)$, \quad & 
$\delta ^+$ & $= {i \over \sqrt 7}\pmatrix 
{\varepsilon - \varepsilon ^6 & \varepsilon ^4 - \varepsilon ^3 & \varepsilon ^2-  
\varepsilon ^5 \cr   
\varepsilon ^4 - \varepsilon ^3 & \varepsilon ^2 - \varepsilon ^5 & \varepsilon -  
\varepsilon ^6 \cr   
\varepsilon ^2- \varepsilon ^5 & \varepsilon  - \varepsilon ^6 & \varepsilon ^4- 
\varepsilon ^3 }$ 
\cr 
& & & & & \cr 
& & & & & \cr 
$\mu ^-$ & $= \pmatrix {1 & 0 & 0 & 0 \cr 
                    0 & 0 & 0 & 1 \cr 
                    0 & 1 & 0 & 0  \cr 
                    0 & 0 & 1 & 0 }, \quad $ & $\nu ^-$ & $= diag\ (1 , \varepsilon ,  
\varepsilon ^ 2, \varepsilon ^4)$, \quad & 
$\delta ^-$ & $= {i \over \sqrt 7}\pmatrix 
{1 & 1 & 1 & 1 \cr 
2 & \varepsilon + \varepsilon ^6 & \varepsilon ^4 +\varepsilon ^3 & \varepsilon ^2 +  
\varepsilon ^5 \cr   
2 & \varepsilon ^4 + \varepsilon ^3 & \varepsilon ^2 +\varepsilon ^5 & \varepsilon + 
\varepsilon ^6 \cr   
2 & \varepsilon ^2+\varepsilon ^5 & \varepsilon +\varepsilon ^6 & \varepsilon ^4+ 
\varepsilon ^3  }$ 
\cr} 
$$

From the character table of $SL_2(\bZ_7)$ one sees that, as a $SL_2(\bZ_7)$-module,  
$V=W'\oplus U'$ and from the above computations one gets concrete realizations of 
$W'$, $U'$, namely  $W'=V^+$ and $U'=V^- $.   
$$ 
\eqalign{ 
S^2W  = &\ T \cr 
S^3W  = &\ L\oplus W' \cr 
S^4W  = &\ I\oplus M_2 \oplus T } 
\quad \quad 
\eqalign{ 
S^2W' = &\ T  \cr 
S^3W'  = &\ L\oplus W  \cr 
S^4W' = &\ I\oplus M_2 \oplus T }  
$$ 
If we denote by  
$$ 
v_1=e_1-e_6  \quad 
v_1=e_4-e_3  \quad
v_1=e_2-e_5  
$$ 
the chosen basis of $W'$, then the only  $SL_2(\bZ_7)$ -- invariant quartic 
is the Klein quartic: 
$$f_{klein}\ =v_1^3v_2+v_2^3v_3+v_3^3v_1\ \ .$$ 
  
$$S^2U'=L\oplus W'\ \ .$$ 
We choose as basis for $L\subset S^2U'$ the following elements: 
$$ 
f_0= u_0^2,\quad
f_1=u_2u_3, \quad
         f_2=u_3u_1, \quad 
         f_3=u_1u_2, \quad
f_4=u_0u_3+u_1^2, \quad
         f_5=u_0u_1+u_2^2, \quad 
         f_6=u_0u_2+u_3^2. 
$$ 
and for $W'$ the elements 
$$ 
v_3=u_0u_3-u_1^2, \quad 
         v_2=u_0u_1-u_2^2, \quad
         v_1=u_0u_2-u_3^2. 
$$ 
Then in the decomposition 
$$S^3V_3=(I\oplus U'\oplus L)V_4$$ 
the elements corresponding to 
$f_je_0$ are given by: 
$$ 
\eqalign {f_0e_0= x_1x_2x_4+x_3x_5x_6-x_0^3\cr 
                    \cr  
                    \cr} 
\quad  
\eqalign{f_1e_0=x_0x_1x_6 \cr 
         f_2e_0=x_0x_2x_5 \cr 
         f_3e_0=x_0x_3x_4 \cr} 
\quad  
\eqalign{f_4e_0=x_2^2x_3+x_5^2x_4 \cr 
         f_5e_0=x_1^2x_5+x_6^2x_2 \cr 
         f_6e_0=x_4^2x_6+x_3^2x_1 \cr 
} 
$$ 
From here one obtains all $f_je_k$ via cyclic permutation, in other words via  
the action of $\sigma$. 
 
\medskip

{\bf References}
\frenchspacing 
\item{[ACGH]} E. Arbarello, M. Cornalba, P.A. Griffith, J. Harris, 
{\it Geometry of Algebraic Curves} Volume I, Grundlehren der 
math. Wiss. {\bf 267}, Springer-Verlag, New York, Berlin, Heidelberg, 
Tokio 1984. 
 
\item{[BE]} D.A. Buchsbaum, D. Eisenbud, {\it Algebra structures for finite free  
resolutions, and some structure theorems for ideals of codimension 3} Am. J.  
Math. 99, (1977) 447-485. 

\item{[BHM]} W. Barth, K. Hulek, R. Moore, {\it Degenerations of  
Horrocks-Mumford surfaces}, Math. Ann. {\bf 277} (1987) 735-755. 
 
\item{[DK]} I. Dolgachev, V. Kanev, {\it Polar covariants of plane cubics 
and quartics}, Adv. Math. {\bf 98} (1993) 216-301.  
 
\item{[E]} D. Eisenbud, {\it Commutative Algebra with a View Toward Algebraic  
Geometry}, Graduate texts in Mathematics {\bf 150}, Springer Verlag ,  
Heidelberg, New York, Berlin, 1994.  
 
\item{[EPW]} D. Eisenbud, S. Popescu, C. Walter, {\it Resolution of 
modules with Symmetry}, in preparation. 
 
\item{[EPS]} G. Ellingsrud, R. Piene, S.A. Str{\o}mme, {\it On the Variety  
of Quadrics Defining Twisted Cubics}, Springer Lect. Notes in Math. 
{\bf 1266}  (1987) 84-96. 
 
\item{[Gr]} V. Gritsenko, {\it Irrationality of moduli spaces of 
polarized abelian surfaces}, in {\it Abelian varieties}, Proceedings of the 
Egloffstein Conference 1993, de Gruyter, Berlin 1995, 63-81. 
 
\item{[GP1]} M. Gross, S. Popescu, {\it Equations of $(1,d)$-polarized 
abelian surfaces }, Math. Ann. {\bf 310} (1998) 333-337. 
 
\item{[GP2]} M. Gross, S. Popescu, {\it Calabi-Yau 3-folds and moduli 
of abelian surfaces I}, in preparation. 
 
\item{[GP3]} M. Gross, S. Popescu, {\it Calabi-Yau 3-folds and moduli 
of abelian surfaces II }, in preparation. 
 
\item{[H]} R. Hartshorne, {\it Algebraic Geometry}, GTM {\bf 52}, 
Springer-Verlag, New York, Heidelberg, Berlin, 1977. 
 
\item{ [HM]} G. Horrocks, D. Mumford, {\it A Vector Bundle on $\bP^4$  
with $15000$ Symmetries}, Topology {\bf 12} (1973) 63-81. 
 
\item{[Isk]} V. A. Iskovskih  Fano 3-folds II, Math. USSR Izv. {\bf 12} (1978) 
469-506. 
 
\item{[K]} F. Klein {\it \"Uber die Transformation siebenter Ordnung der  
elliptischen Funktionen} Math, Annalen, {\bf 14} (1878) , or the same title in 
{\it Gesammelte Mathematische Abhandlungen} {\bf 3.} Band, Verlag von Julius 
Springer, Berlin, 1923. 
 
\item{[LB]} H. Lange and Ch. Birkenhake, 
{\it Complex abelian varieties}, Grundlehren der math. Wiss. 
{\bf 302} Springer-Verlag, New York, Heidelberg, Berlin, 1992. 
 
\item{[Man86]} N. Manolache,  {\it On the Normal Bundle to Abelian Surfaces  
Embedded in $\bP^4(\bC)$}, Manuscripta math. {\bf 55} (1986) 111-119.  
 
\item{[Man89]} N. Manolache, {\it Syzygies of Abelian Surfaces  
Embedded in $\bP^4(\bC)$}, J. reine angew. Math. {\bf 384} (1988) 180-191 , 
{\it The Equations of the Abelian Surfaces Embedded in $\bP^4(\bC)$}, J. reine 
angew. Math. {\bf  394} (1989) 196-202. 
 
\item{[Muk89]} S. Mukai, {\it Biregular Classification of Fano 3-folds  
and Fano manifolds  of coindex 3}, Proc. Natl. Acad. Sci. USA {\bf 86} (1989)  
3000-3002.  
 
\item{[Muk92]} S. Mukai, {\it Fano 3-folds}, Sel. Pap. Conf. Proj. Var.,  
Trieste/Italy 1989, and Vector Bundles and Special Proj. Embeddings,  
Bergen/Norway 1989, Lond. Math. Soc. Lect. Notes Ser. {\bf 179} (1992) 255-263. 
 
\item{[Mum66]}  D. Mumford, {\it On the Equations Defining Abelian  
Varieties. I }, Invent. math. {\bf 1} (1966) 287--354.  
 
\item{[Mum70]} D. Mumford,  {\it Abelian Varieties}, Oxford University  
Press 1966. 
 
\item{[P]} C. Plegge, {\it Gleichungen f\"ur Modulkurven}, 
Dissertation Bayreuth, 1997.  
 
\item{[PS]} R. Piene, M. Schlessinger, {\it On the Hilbert Scheme  
Compactification of the Space of Twisted Cubics}, Amer. J. Math. {\bf 107}   
(1985) 761-774. 
 
\item{[Ros]} J. Rosanes, {\it Ueber ein Princip der Zuordnung algebraische  
Formen}, J. reine und angew. Math. {\bf 76} (1873) 312-330  
 
\item{[R]} I. Reider, {\it Vector bundles of rank 2 and linear systems on  
algebraic surfaces}, Ann. Math. {\bf 127} (1988) 309-316. 
 
\item{[Schr]} F.-O. Schreyer, {\it Geometry and Algebra of Prime Fano  
3-folds of genus 12}, in preparation.  
 
\item{[Sco1]} G. Scorza, {\it Un Nuovo Teorema sopra le Quartiche  
Piane Generali}, Math. Ann. {\bf 52} (1889). 
 
\item{[Sco2]} G. Scorza, {\it Sopra la Teoria delle Figure Polari 
delle Curve Piane del $4^o$ Ordine}, Ann. di Mat. (3) {\bf 2} (1899) 155-202.  
 
\item{[W]} C. Walter, {\it Pfaffian subschemes}, J. Alg. Geom. {\bf 5} 
(1996) 671-704. 
           
\medskip 
\noindent Authors' addresses: 
\medskip 
 Nicolae Manolache \par
{\it Institute of Mathematics \par
Romanian Academy \par
P. O. Box 1--764, RO-70700 \par
Bucharest, Romania \par}
e-mail: manole@stoilow.imar.ro \par
\noindent current address: \par
{\it Fachbereich Mathematik \par
Universit\"at Oldenburg \par
Pf. 2503          \par
D-26111 Oldenburg \par
Germany} \par
e-mail: manolache@math.uni-oldenburg.de \par
 
\medskip 
 Frank-Olaf Schreyer \par 
{\it Mathematisches Institut \par  
Universit\"at Bayreuth \par 
D-95440 Bayreuth \par 
Germany \par} 
e-mail: frank.schreyer@uni-bayreuth.de 
   
\bye